\documentclass[12pt]{article}

\usepackage{latexsym}
\usepackage{amsfonts}
\usepackage{amsmath}
\usepackage{amsthm}
\usepackage{epsfig}

\newcommand {\ve} {{\mbox{ and }}}
\newcommand{\hamil} {{\cal H}}
\newcommand{\cov} {\mbox{\textup{cov}}}
\newcommand{\prls} {Peierls}
\newcommand{\nsim} {{\not\sim}}
\newcommand{\I} {{\mbox{\bf I}}}

\newcounter{AbcT}

\newtheorem {Theorem}    {Theorem}[section]
\newtheorem {Definition} {Definition} [section]
\newtheorem {Problem}    {Problem}

\newtheorem {Lemma}      [Theorem]    {Lemma}
\newtheorem {Corollary}  [Theorem]    {Corollary}
\newtheorem {Proposition}[Theorem]    {Proposition}
\newtheorem {Claim}      [Theorem]    {Claim}

\newtheorem {Conjecture}[Theorem]     {Conjecture}

\newcommand{\Z}{{\bf{Z}}}
\newcommand{\ignore}[1]{}

\def\sT{n}
\def\sTrb{n_r}
\def\sGr{n_r}
\def\Vol{n}

\def\eps{\epsilon}
\def\E{{\bf{E}}}
\def\P{{\bf{P}}}
\def\Var{{\bf{Var}}}

\def\R{\hbox{I\kern-.2em\hbox{R}}}
\def\A{{\cal{A}}}

\def\F{{\cal{F}}}

\def\|{ \Vert  }
\def\v0{{\bf 0}}
\def\one{{\bf 1}}
\def\0{\hat{0}}
\def\1{\hat{1}}
\def\al{\alpha}

\def\path{{\tt path}}

\def\lam{\lambda}

\def\T{{\bf T}}
\def\phi{\varphi}

\def\be{\begin{equation}}
\def\ee{\end{equation}}

\def\one{{\bf 1}}
\def\eps{\epsilon}
\def\lam{\lambda}
\def\cE{{\xi}}
\def\dir{{\cal E}}

\def\Bin{{\sf Bin}}
\def\L{{\mathcal L}}

\begin{document}
\title{Glauber Dynamics on Trees \\ and Hyperbolic Graphs}
\author{
  Noam Berger \thanks{Research supported by Microsoft graduate fellowship.}
   \\ University of California, Berkeley \\
  noam@stat.berkeley.edu \and
  Claire Kenyon \\ LRI, UMR CNRS \\
  Universit\'e Paris-Sud, France \\kenyon@lix.polytechnique.fr  \and
  Elchanan Mossel \thanks{Supported by a visiting position at INRIA
    and a PostDoc at Microsoft research.}
  \\ University of California, Berkeley \\
  mossel@stat.berkeley.edu \and
  Yuval Peres\thanks{ Research supported by NSF Grants
DMS-0104073, CCR-0121555 and a Miller Professorship at UC Berkeley.}
 \\ University of California, Berkeley \\
  peres@stat.berkeley.edu
}

\maketitle
\thispagestyle{empty}

\begin{abstract}
We study continuous time Glauber dynamics for random configurations with
local constraints (e.g. proper coloring, Ising and Potts models) on
finite graphs with $n$ vertices and of bounded degree.
We show that the relaxation time
 (defined as the reciprocal of the spectral gap
$|\lambda_1-\lambda_2|$) for the dynamics on trees and on planar
hyperbolic graphs, is polynomial in $n$. For these hyperbolic graphs,
this yields a general polynomial sampling algorithm for random configurations.
We then show that for general graphs, 
if the relaxation time $\tau_2$ satisfies $\tau_2=O(1)$,
then the correlation coefficient, and the mutual information, between
any local function (which depends only on the configuration in a fixed
window) and the boundary conditions, decays exponentially in the distance
between the window and the boundary.
For the Ising model on a regular tree, this condition is sharp.

\end{abstract}
\section{Introduction}

\subsubsection*{Context}
In recent years, Glauber dynamics on the lattice $\Z^d$ was
extensively studied. A good account can be found in \cite{Ma}.
In this work, we study this dynamics on other graphs.

The main goal of our work is to determine
which geometric properties of the underlying graph
 are most relevant to the mixing rate
of the Glauber dynamics on particle systems.

To define a general particle system~\cite{Li}
on an undirected graph $G=(V,E)$,
define a configuration as an element $\sigma$ of $A^{V}$ where
$A$ is some finite set, and to each edge $(v,w) \in E$,
associate a weight function
 $\al_{vw} : A\times A \to \R_{+}$.
The Gibbs distribution assigns every configuration $\sigma$ probability
proportional to $\prod_{ \{ v,w\} \in E}
\al_{vw}(\sigma_v,\sigma_w)$.
The Ising model
(for which $\al_{vw}(\sigma_v,\sigma_w)=e^{\beta\sigma_v\sigma_w}$)
 and the Potts model are examples of such systems;
so is the coloring model (for which
$\al_{vw}=\one_{\sigma_v\neq\sigma_w}$)

On a finite graph, the Heat-Bath Glauber dynamics is a continuous
time Markov chain with the generator
\begin{equation}\label{generator}
(\L(f))(\sigma)=\sum\limits_{v\in V}
\left(\sum\limits_{a\in A} K[\sigma \to \sigma_v^a]
\left(f(\sigma_v^a)-f(\sigma)\right) \right),
\end{equation}
where $\sigma_v^a$ is the configuration s.t.
\begin{equation*}
\sigma_v^a(w) = \left\{
\begin{array}{ll}
a & \mbox{ if }w=v \\
\sigma(w) & \mbox{ if }w\neq v
\end{array}
\right.
\end{equation*}
and

\begin{equation*}
K[\sigma \to \sigma_v^a]=\frac
{\prod\limits_{w:(w,v)\in E}
\al_{vw}(a,\sigma_w)}
{\sum\limits_{a^\prime\in A}\left({\prod\limits_{w:(w,v)\in E}
\al_{vw}(a^\prime,\sigma_w)}\right)}.
\end{equation*}

It is easy to check that this dynamics is reversible with respect to the
Gibbs measure.
An equivalent representation for the Glauber dynamics, known as the
{\em Graphical representation}, is the following:
Each vertex has a rate $1$ Poisson clock attached to it. These Poisson
clocks are independent of each other.
Assume that the clock at $v$ rang at time $t$ and that just before
time $t$ the configuration was $\sigma$.
Then at time $t$ we replace $\sigma(v)$ by a random spin $\sigma'(v)$
chosen according to the Gibbs distribution conditional on the
rest of the configuration:
\begin{equation*}
\frac{\P[ \sigma'(v) = i \mid \sigma]}{\P[\sigma'(v) = j \mid \sigma]}
= \prod_{w: \{ v,w\} \in E}
{\al_{vw}(i,\sigma(w))\over \al_{vw}(j,\sigma(w))} .
\end{equation*}

We are interested in the rate of convergence of the Glauber dynamics
to the stationary distribution. Note that this process mixes $n=|V|$
times faster than
the corresponding discrete time process, simply because it performs (on
average) $n$ operations per time unit while the discrete time process
performs one operation per time unit.

\begin{figure} [hptp] \label{5gon}
\epsfxsize=5cm
\hfill{\epsfbox{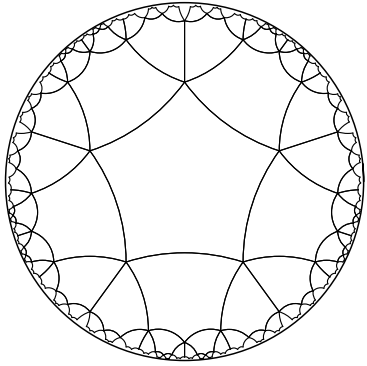}\hfill}
\caption{A ball in hyperbolic tiling}
\end{figure}

In section \ref{subsec:cut-width},
we describe a connection between the geometry of a graph and the mixing time
of Glauber dynamics on it.
In particular, we show that
for balls in hyperbolic tilings, the Glauber dynamics for the Ising model,
the Potts model and proper coloring with $\Delta+2$ colors
(where $\Delta$ is the maximal degree),
have mixing time polynomial in the volume.
An example of such a hyperbolic  graph
can be obtained from the binary tree by adding horizontal edges across levels;
another example is given in Figure \ref{5gon}.

In sections \ref{subsec:recur}-\ref{sec:hightemp}
we study Glauber dynamics for the Ising model
on regular trees.
For these trees
we show that the mixing time is polynomial at all temperatures, and
we characterize
the range of temperatures for which the spectral gap
is bounded away from zero.
Thus, the notion that the two sides of the phase transition
(high versus low temperatures) should correspond to polynomial versus
super-polynomial  mixing times for the associated dynamics,
fails for the Ising model on trees: here the two sides of the high/intermediate
versus low temperature phase transition just correspond to uniformly bounded versus
{\em unbounded} inverse spectral gap. We
also exhibit another surprising phenomenon:
On infinite regular trees, there is a range
of temperatures in which the inverse spectral gap is bounded, even though
there are many different Gibbs measures.

In section \ref{sec:general} of the paper we go beyond trees and
hyperbolic graphs and study Glauber dynamics for
families of finite graphs of bounded degree.
We show that if the inverse spectral
gap of the Glauber dynamics on the ball centered at $\rho$
stays bounded as the ball grows,
then the correlation between the
state of a vertex $\rho$ and the states of
vertices at distance $r$ from $\rho$, must
 decay exponentially in $r$.

\subsubsection*{Setup}
\noindent{\bf The graphs.}
Let $G=(V,E)$ be an infinite graph with maximal degree $\Delta$.
Let $\rho$ be a distinguished vertex and denote by $G_r = (V_r,E_r)$ the
induced graph
on $V_r = \{v \in V : {\rm dist}(\rho,v) \leq r\}$.
Let $\sGr$ be the number of vertices in $G_r$.
At some parts of the paper we will
focus on the case where $G=T=(V,E)$ is the infinite
$b$-ary tree.
In these cases, $T_r^{(b)}=(V_r,E_r)$ will denote the $r$-level $b$-ary tree.

\noindent{\bf The Ising model.}
In the Ising model on $G_r$ at inverse temperature $\beta$, every
configuration $\sigma \in \{ -1,1 \}^{V_r}$ is assigned probability
$$
\mu[\sigma]=Z(\beta)^{-1} \exp\Big(\beta \sum_{\{v,w\}\in E_r} \sigma(v)
\sigma(w) \Big)
$$
where    
$Z(\beta)$ is a normalizing constant.
When $G_r = T_r^{(b)}$, this measure has the following equivalent
definition~\cite{EKPS}:
Fix $\epsilon=(1+e^{2 \beta})^{-1}$.
Pick a  random spin $\pm 1$
uniformly for the root of the tree.
Scan the tree top-down, assigning vertex $v$ a spin equal to
the spin of its parent with probability $1-\epsilon$ and opposite with
probability $\epsilon$.

The Heat-Bath Glauber dynamics for the Ising model
chooses the new spin  $\sigma'(v)$ in such a way that:
\begin{equation*}
\frac{\P[ \sigma'(v) = +1 \mid \sigma]}{\P[\sigma'(v) = -1 \mid \sigma]}
= \exp\Big( 2 \beta \sum_{w: \, \{ w, v\} \in E_r}\sigma(w)\Big) \,.
\end{equation*}
See \cite{Li} or \cite{Ma} for more background.

\vspace{0.5in}
\noindent{\bf Mixing times.}
\begin{Definition}
For a reversible continuous time Markov chain,
let $0 = \lam_1 \leq \lam_2 \leq \ldots \leq \lam_k$
be the eigenvalues of $-\L$ where $\L$ is the generator.
The {\bf spectral gap} of the chain is
defined
as $\lam_2$, and the {\bf relaxation time},
$\tau_2$, is defined as
the inverse of the spectral gap.
\end{Definition}
Note that the corresponding discrete time Glauber dynamics has
transition matrix $M=\I+{\frac{1}{n}}\L$, where $n$ is the number 
of vertices. Moreover, the eigenvalues of $M$ are
$1,1-{\frac{\lam_2}{n}},1-{\frac{\lam_3}{n}},\ldots$ and
therefore the spectral gap of the discrete dynamics is the spectral
gap of the continuous dynamics divided by $n$.
\begin{Definition}
For measures $\mu$ and $\nu$ on the same discrete space,
the {\bf total-variation distance}, $d_V(\mu,\nu)$, between $\mu$ and $\nu$
is defined
as
\[
d_V(\mu,\nu) = \frac{1}{2} \sum_{x} |\mu({x}) - \nu({x})| \,.
\]
\end{Definition}
\begin{Definition}
Consider an ergodic Markov chain $\{X_t\}$
with stationary distribution $\pi$ on a finite state
space. Denote by $\P^t_x$ the law of $X_t$ given $X_0=x$.
The {\bf mixing time} of the chain, $\tau_1$, is defined as
\[
\tau_1 = \inf \{ t : \sup_{x,y} d_V(\P^t_x,\P^t_y) \leq e^{-1}\}.
\]
\end{Definition}
\noindent For $t \ge \ln (1/\epsilon )\tau_1$, we have
$$\sup_{x} d_V(\P^t_x,\pi) \leq \sup_{x,y} d_V(\P^t_x,\P^t_y) \leq
\epsilon .$$

Using $\tau_2$ one can bound the mixing time $\tau_1$, since every
reversible Markov
chain with stationary distribution $\pi$ satisfies (see, e.g., \cite{Al}),
\begin{equation} \label{eq:t1_and_t2}
\tau_2 \leq \tau_1 \leq \tau_2
\left( 1 + \log \left( (\min_{\sigma} \pi(\sigma))^{-1} \right)
\right).
\end{equation}
For the Markov chains studied in this paper, this gives
$\tau_2 \leq \tau_1 \leq O(n)\tau_2$.

\vspace{0.25in}
\noindent{\bf Cut-Width and relaxation time.}
\begin{Definition}
The {\bf cut-width} $\cE(G)$ of a graph $G$ is the
smallest integer such that
there exists a labeling  $v_1,\ldots,v_n$ of the vertices
such that for all $1 \leq k \leq n$,
the number of edges from $\{v_1,\ldots,v_k\}$ to
$\{v_{k+1},\ldots,v_n\}$, is at most  $\cE(G)$.
\end{Definition}
\begin{Remark}
The {\em vertex-separation} of a graph $G$
is defined analogously to the cut-width in terms
of vertices among $\{v_1,\ldots,v_k\}$ that are adjacent to
$\{v_{k+1},\ldots,v_n\}$. In \cite{Ki} it is shown that the vertex-separation of
$G$ equals its {\em path-width}, see \cite{RS}.
In \cite{first} the cut-width was called the {\em exposure}.
\end{Remark}

Generalizing an argument in \cite[Theorem 6.4]{Ma} for
$\Z^d$, (see also \cite{JS1}), we prove:
\begin{Proposition} \label{prop:cut-width}
Let $G$ be a finite graph with $n$ vertices and
maximal degree $\Delta$.

\begin{enumerate}
\item\label{prop:exp:ising}
Consider the Ising model on $G$.
The relaxation time of the Glauber dynamics is at most
$n e^{(4 \cE(G) + 2 \Delta) \beta}$.
\item\label{prop:exp:color}
Consider the coloring model on $G$. If the number of colors $q$
satisfies $q \geq \Delta + 2$, then the relaxation time of the Glauber
dynamics is at most
$(\Delta + 1) n (q-1)^{\cE(G) + 1}$.
\end{enumerate}
\end{Proposition}
\noindent
Analogous results hold for the independent set and hard core models.

\noindent {\bf Cut-Width and long-range correlations for hyperbolic
graphs.} The usefulness of Proposition \ref{prop:cut-width} comes
about when we bound the relaxation time of certain graphs by
estimating their cut-width. The following proposition bounds the
cut-width of balls in hyperbolic tilings of the plane. Recall that
the {\bf Cheeger constant} of an infinite graph $G$ is
\begin{equation}\label{eq:def_cheeg}
c(G) =
\inf\left\{\left.
\frac{|\partial A|}{|A|}\right|
A\subseteq G ; 0 < |A| < \infty
\right\} ,
\end{equation}
where $\partial A$ is the set of vertices of $A$ which have
neighbors in $G\setminus A$.

\begin{Proposition}\label{prop:hyptylexp}
For every $c > 0$ and $\Delta < \infty$,
there exists a constant $C = C(c,\Delta)$ such that if $G$ is an infinite planar graph with 
\begin{itemize}
\item
Cheeger constant at least $c$, 
\item
maximum degree bounded by $\Delta$ and 
\item
for every $r$ no cycle from $G_r$ separates two
vertices of $G\setminus G_r$,
\end{itemize} 
then $\cE(G_r) \leq C
\log n_r$ for all $r$, where $n_r$ is the number of vertices of $G_r$. 
\end{Proposition}
Combining this with Proposition \ref{prop:cut-width} we get that the
Glauber dynamics for the Ising models on balls in the hyperbolic
tiling has relaxation time polynomial in the volume for every
temperature. On the other hand, we have the following proposition:
\begin{Proposition}\label{prop:corrhyp}
Let $G$ be a planar graph with bounded degrees, bounded co-degrees and a
positive Cheeger constant. Then there exist $\beta' < \infty$ and
$\delta > 0$ such that for
all $r$, all $\beta > \beta'$, and all vertices $u,v$ in $G_r$,
the Ising model on $G_r$ satisfies that $\E[\sigma_u \sigma_v] \geq \delta$.
In other words, at low enough temperature there are
long-range correlations.
\end{Proposition}
This shows that for the Ising model on balls of hyperbolic tilings at
very low temperature, there are long-range correlations coexisting with
polynomial time mixing. 
While there is no characterization of all planar graphs with positive Cheeger constants, an important family of such graphs is  
Cayley graphs of nonelementary fuchsian groups which are nonamenable, and hence have a positive Cheeger constant, see \cite{pater,magnus,katok} for background and \cite{HJR} for some
explicit estimates.

\smallskip

\noindent{\bf Relaxation time for the Ising model on the tree.}
The Ising model on the
$b$-ary tree has three different regimes,
see \cite{BRZ, EKPS}.
In the high temperature regime, where $1 - 2 \eps < 1/b$, there is
a unique Gibbs measure on the infinite tree, and the expected value of the spin at the
root $\sigma_{\rho}$ given any
boundary conditions $\sigma_{\partial T_r^{(b)}}$ decays
exponentially in $r$.
In the intermediate regime, where $1/b < 1 - 2 \eps < 1/\sqrt{b}$,
the exponential decay described above still holds
for typical boundary conditions, but not for certain exceptional
boundary conditions, such as the all $+$ boundary;
consequently, there are
infinitely many Gibbs measures on the infinite tree.
In the low temperature regime, where $1 - 2 \eps > 1/\sqrt{b}$, typical
boundary conditions impose bias on the expected value of the spin at the root $\sigma_{\rho}$.

\begin{Theorem} \label{thm:tree}
Consider the Ising model on the $b$-ary tree $T_r= T_r^{(b)}$ with $r$ levels.
Let $\epsilon=(1+e^{2 \beta})^{-1}$.
The relaxation time $\tau_2$
for Glauber dynamics on $T_r^{(b)}$ can be bounded as follows:
\begin{enumerate}
\item \label{treeupperbound}
The relaxation time is polynomial at all temperatures:
$\tau_2=n_r^{O(\log(1/\epsilon ))}$. Furthermore, the limit
\[
\lim\limits_{r\to\infty}\frac{\log(\tau_2(T_r^{(b)},\beta))}{\log(n_r)}
\]
exists.
\item {Low temperature regime:}
        \begin{enumerate}
        \item \label{treesuperlinear}
If $1-2\epsilon\geq 1/\sqrt{b}$ then $\sup_r \tau_2 (T_r) =\infty$.
In fact, $\tau_2(T_r) =\Omega(n_r^{\log_b(b(1-2\epsilon )^2)})$ when 
$1 - 2\epsilon > 1/\sqrt{b}$ and $\tau_2(T_r) = \Omega(\log n_r)$ when
$1 - 2\epsilon = 1/\sqrt{b}$.  
        \item \label{treeinfinitedegree}
Moreover, the degree of $\tau_2$  tends to infinity
as $\eps$ tends to zero:
 $\tau_2(T_t) =n_r^{\Omega(\log (1/\epsilon ))}$.
        \end{enumerate}
\item \label{treelinear} {Intermediate and high temperature regimes:} \\
If $ 1-2\epsilon < 1/\sqrt{b} $ then
the relaxation time is uniformly bounded:
$\tau_2 = O(1)$. Furthermore, this result holds for every external field
$\{H(v)\}_{v\in T_r}$.
\end{enumerate}
\end{Theorem}

In particular we obtain from Equation (\ref{eq:t1_and_t2})
that in the low temperature region
$\tau_1 = \sTrb^{\Theta(\beta)}$, and in the intermediate and high
temperature regions $\tau_1 = O(\sTrb)$.

A recent work by Peres and Winkler \cite{PW03} compares
the mixing times of single site and block dynamics
for the heat-bath Glauber dynamics for the Ising model.

They show that if the blocks are of bounded volume, then
the same mixing time up to constants is obtained for the single site
and block dynamics.

Combining these results with the path coupling argument of Section \ref{sec:hightemp},
it follows that
$\tau_1 = O(\log \sTrb)$ in the intermediate and high
temperature regions.

We emphasize that Theorem \ref{thm:tree} implies that in the
intermediate region
$1/2 < 1 - 2 \eps < 1/\sqrt{b}$,
the relaxation time is bounded by a constant,
yet, in the infinite volume
there are infinitely many Gibbs measures.
This Theorem is perhaps easiest to appreciate when compared to other
results on the Gibbs distribution for the Ising model
on binary trees, summarized in Table~\ref{tableau}.

\begin{table}\label{tableau}
\begin{center}
{\small
\begin{tabular}{|l|c|l|l|c|}\hline
Temp. & $1-2\epsilon$ & $\sigma_{\rho}|\sigma_{\partial T}\equiv +   $ &
$I(\sigma_{\rho},\sigma_{\partial T})$
               &  $\tau_2$ \\ \hline\hline
high & $< 1/2$ & unbiased & $\to 0$ & $O(1)$\\ \hline
med. & $\in (\frac{1}{2}, \frac{1}{\sqrt{2}})$ & biased
    & $\to 0$ & $O(1)$\\ \hline
low & $> \frac{1}{\sqrt{2}}$ & biased & inf$>0$ &
    $n^{\Omega(1)}$ \\ \hline
freeze & $1-o(1)$ & biased & $1 - o(1)$ &
    $n^{\Theta(\beta)}$ \\ \hline
\end{tabular} }
\end{center}
\caption{The Ising model on binary trees. Here the root is denoted $\rho$,
and the vertices at distance $r$ from the root are denoted $\partial T$.}
\end{table}
The proof of the low temperature result is quite general and applies
to  other models with ``soft'' constraints,
such as Potts models on the tree .

\noindent{\bf Spectral gap and correlations.}
At infinite temperature, where distinct vertices are independent, the
 Glauber dynamics on a graph of $n$ vertices
reduces to an (accelerated by a factor of $n$) random walk on a discrete $n$-dimensional cube,
where it is well known that the relaxation time is $\Theta(1)$.
Our next result shows that at any temperature where
such fast relaxation takes place, a strong form of independence holds.
This is well known in $\Z^d$, see~\cite{Ma}, but our formulation is valid
for any graph of bounded degree.
\begin{Theorem} \label{thm:general}
Denote by $\sigma_r$ the configuration on all vertices
at distance
$r$ from $\rho$.
If $G$ has bounded degree and the relaxation time
of the Glauber dynamics satisfies $\tau_2(G_r) = O(1)$, then
the Gibbs distribution on $G_r$ has the following property.
For any fixed finite set of vertices $A$, there exists $c_A>0$ such
that for $r$ large enough
\begin{equation} \label{eq:cor_decay}
{\rm Cov}[f,g] \le  e^{-c_A r} \sqrt{{\rm Var}(f){\rm Var}(g)} \, ,
\end{equation}
provided that $f(\sigma)$ depends only on $\sigma_A$ and $g(\sigma)$ depends
only on $\sigma_r$. Equivalently, there exists $c_A'>0$ such that
\begin{equation} \label{eq:mut_inf_decay}
I[\sigma_A,\sigma_r] \le e^{-c_A' r} \, ,
\end{equation}
where $I$ denotes mutual information (see~\cite{CT}.)
\end{Theorem}
This theorem holds in a very general setting which includes Potts models,
random colorings, and other local-interaction models.

Our proof of Theorem \ref{thm:general} uses ``disagreement
percolation'' and a  coupling argument
exploited by van den Berg, see \cite{vB}, to establish uniqueness of
Gibbs measures in $\Z^d$; according to F. Martinelli (personal communication)
 this kind of argument is originally  due to B. Zegarlinski.
Note however, that Theorem \ref{thm:general} holds also when there
are multiple Gibbs measures --
as the case of the Ising model in the intermediate regime demonstrates.
Moreover, combining Theorem \ref{thm:general} and Theorem \ref{thm:tree}, one
infers that for $1 - 2 \eps < 1/\sqrt{b}$, we have
$\lim_{r \to \infty} I[\sigma_0,\sigma_r] = 0$.
This yields another proof
of this fact which was proven before in \cite{BRZ,Io,EKPS}.

\subsubsection*{Plan of the paper}
In section~\ref{sec:alltemp} we prove Proposition~\ref{prop:cut-width}
via a canonical path argument, and give the resulting polynomial time
upper bound of
Theorem~\ref{thm:tree} part \ref{treeupperbound}.
We also present a more elementary proof of the upper bound
on the relaxation time for the tree, which gives sharper exponents and
the existence of a limiting exponent;
this
proof uses Martinelli's block dynamics to
show sub-additivity.
In section~\ref{sec:lowerbounds} we sketch a proof of
Theorem~\ref{thm:tree} part~\ref{treesuperlinear} and present a proof
of Theorem~\ref{thm:tree} part \ref{treeinfinitedegree}.
These  lower bounds are
 obtained by finding a low conductance ``cut'' of the configuration
space, using global majority of the boundary spins
for the former result,
and recursive majority for the latter result.
In section~\ref{sec:hightemp}
we establish the high temperature result, using comparison to
block dynamics which are analyzed via path-coupling.
Finally, in section~\ref{sec:general} we prove Theorem~\ref{thm:general} by a Peierls
argument controlling ``paths of disagreement''
between two coupled dynamics.

{\bf Remark:}
Most of the results proved here were presented
(along with proof sketches) in the extended abstract \cite{first}.
However, the proofs of our results for hyperbolic graphs
(see Section 2.2), which involve some interesting geometry,
 were not even sketched there. Also, the general polynomial
upper bound for trees that we establish in Section 2.3 
is a substantial improvement 
on the results of \cite{first}, since it only assumes the dynamics is
  ergodic and allows for arbitrary hard-core constraints.

\section{Polynomial Upper Bounds}\label{sec:alltemp}
\subsection{Cut-Width and mixing time} \label{subsec:cut-width}

We begin by showing how part
\ref{prop:exp:ising} of Proposition \ref{prop:cut-width} implies the upper
bound in part \ref{treeupperbound} of Theorem~\ref{thm:tree}.

\begin{Lemma}\label{firstdepthsearch}
Let $T_r^{(b)}$ be the $b$-ary tree with $r$ levels. Then,
$\cE(T_r^{(b)})<(b-1)r+1.$
\end{Lemma}
\begin{proof}
Order the vertices using the {\em Depth first search left to right
  order }, i.e., use the
following labeling for the vertices: The root is labeled
$\langle 0,0,\ldots,0\rangle$. The children of the root are labeled $\langle 1,0,\ldots,0\rangle$
through $\langle b,0,\ldots,0\rangle$, and so on, so that the children of
$\langle a_1,a_2,\ldots,a_k,0,\ldots,0\rangle$ are $\langle a_1,a_2,\ldots,a_k,1,\ldots,0\rangle$
through $\langle a_1,a_2,\ldots,a_k,b,\ldots,0\rangle$. Then order the vertices
lexicographically. Note that in the lexicographic ordering, a vertex always appears before its children.
When we enumerated all vertices up to
$\langle a_1,a_2,\ldots,a_r\rangle$, the only vertices that were enumerated but whose
children were not enumerated are among the set of at most $r$ vertices
\[
\left\{
\langle 0,0,\ldots,0\rangle,\langle a_1,0,\ldots,0\rangle,\langle a_1,a_2,\ldots,0\rangle,\ldots,\langle a_1,a_2,\ldots,a_r\rangle
\right\}.
\]
Each of these vertices has at most $b$ children,
and for all but $\langle a_1,a_2,\ldots,a_r\rangle$ at least one child has
already been enumerated. Therefore,
$$\cE(T_r^{(b)})<(b-1)r+1.$$
\end{proof}

\begin{Corollary}\label{cor:thm:tree}
\begin{enumerate}
\item
The relaxation time of the Glauber dynamics for the Ising model on $T_r^{(b)}$ is at most
\[
C(\eps) \sTrb^{1 + 4(b-1) \log_b\frac{1-\eps}{\eps}}=n_r^{O(\log
  (1/\epsilon ))}.
\]
\item
The relaxation time of the Glauber dynamics for the coloring on
$T_r^{(b)}$ with $q>b+2$ colors is at most
\[
(b+1)\sTrb^{1 + 2(b-1) \log_b(q)}
\]
\end{enumerate}
\end{Corollary}
\begin{proof}
The Corollary follows from Lemma \ref{firstdepthsearch} and
Proposition \ref{prop:cut-width}.
\end{proof}
The upper bound in part~\ref{treeupperbound} of Theorem~\ref{thm:tree}
follows immediately.

\begin{proof}[Proof of part \ref{prop:exp:ising} of Proposition
  \ref{prop:cut-width}]
The proof follows the lines
of the proof given in \cite[Theorem 6.4]{Ma} for
the Ising model in $\Z^d$, (see also \cite{JS1}).

Let $\Gamma$ be the graph corresponding to the transitions of the
Glauber dynamics  on the graph $G$. Between any two configurations
$\sigma$ and $\eta$, we define a ``canonical path" $\gamma(\sigma
,\eta )$ as follows. Fix an order $<$ on the vertices of $G$ which
achieves the cut-width. Consider the vertices $v_1 < v_2 < \ldots$
at which $\sigma_v \neq \eta_v$.

We define the $k$-th configuration $\sigma^{(k)}$ on the path
$\gamma(\sigma,\eta)$ by giving spin $\sigma_v$ to every
labeled vertex $v\leq v_k$, spin $\eta_v$ to every labeled vertex $v>v_k$,
and spin $\sigma_v=\eta_v$ for every unlabeled vertex $v$. Note
that $\sigma^{(0)}=\eta$ and  $\sigma^{(d(\sigma ,\eta
))}=\sigma$. Since $\sigma^{(k-1)}$ and $\sigma^{(k)}$ are
identical except for the spin of vertex $v_k$, they are adjacent
in $\Gamma$. This defines $\gamma(\sigma ,\eta )$ (see
Figure~\ref{fig:can-path}).
\begin{figure}[b]
\center{\epsfig{figure=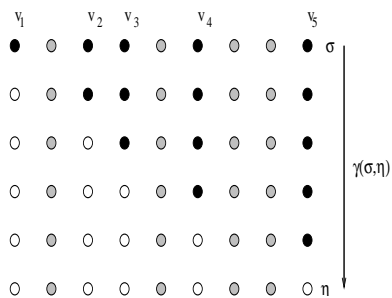, height=1.5in, width=2in}}
\caption{The canonical path from $\sigma$ to $\eta$. The vertices
on which $\sigma$ and $\eta$ agree are marked in grey; the other
vertices are  colored black if their spin is chosen according to
$\sigma$ and white if their spin is chosen according to $\eta$.
}\label{fig:can-path}
\end{figure}
Note that for every $k$, there are at
most $\cE(G)$ pairs  of
adjacent vertices  $(v_i,v_j)$ such that $i\leq k<j$, hence any
configuration on the canonical path between $\sigma$ and $\eta$ will have at
most $\cE(G)$  edges between spins copied from
$\sigma$ and spins copied from $\eta$.

{\noindent \bf Using canonical paths to bound the mixing rate.}
For each directed edge ${\bf e}=(\omega,\zeta)$ on the
configuration graph $\Gamma$, we say that ${\bf e}\in
\gamma(\sigma ,\eta)$ if $\omega$ and $\zeta$ are adjacent
configurations in $\gamma(\sigma ,\eta).$ Let
$$\rho = \sup_{{\bf e}}
\sum_{\sigma , \eta : \, {\bf e}\in \gamma(\sigma ,\eta)}
{\frac{\mu[\sigma]\mu[\eta]}{Q({\bf e})}},$$
where
$\mu$ is the stationary measure (i.e. the Gibbs distribution),
and for any two adjacent configurations $\omega$ and $\zeta$,
$Q({\bf e})=Q((\omega,\zeta))=\mu[\omega]K[\omega \to \zeta]$.
If $L$ is the maximal length of a canonical path, then
by the argument in \cite{JS1,Ma}, the relaxation time of the Markov
chain is at most
\begin{equation} \label{eq:canon_path}
\tau_2 \leq  L \rho.
\end{equation}
Since $L \leq n$, it follows that $\tau_2 \leq n\rho$,
thus it only remains to prove an upper bound on $\rho$.

{\bf Analysis of the canonical path}. For each directed edge ${\bf
e}$ in $\Gamma$, we define an injection from canonical paths going
through ${\bf e}$ in the specified direction, to configurations on
$G$. To a canonical path  $\gamma(\sigma ,\eta)$ going through
${\bf e}$, such that ${\bf e}=(\sigma^{(k-1)},\sigma^{(k)})$, we
associate the configuration $\phi$ which has spin $\eta_{v_i}$ for
every $v_i$ s.t. $i \geq k$ and spin $\sigma_{v_i}$ for every
$v_i$ s.t.  $i < k$. To verify that this  is an injection, note
that one can reconstruct $\sigma$ and $\eta$ by first identifying
the unique $k_0$ s.t. $\omega$ and $\zeta$ differ on $v_{k_0}$ and
then taking (as in Figure~\ref{fig:injection})
\begin{figure}
\center{\epsfig{figure=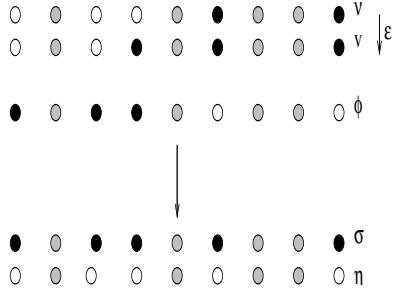, height=1.5in, width=2in}}
\caption{The injection from $({\bf e},\phi )$ to $(\sigma , \eta )
$. The vertices on which both endpoints of ${\bf e}$ and $\phi$
agree are marked in grey; the other vertices are colored black if
they precede $v_{k_0}$ and their spin is chosen according to
$\phi$, or if they are preceded by $v_{k_0}$ and their spin is
chosen according to the endpoints of ${\bf e}$; and are colored
white otherwise.}\label{fig:injection}
\end{figure}
\[
\sigma_{v_k} = \left\{
\begin{array}{ll}
\omega_{v_k} & \omega_{v_k} = \phi_{v_k} \\
\omega_{v_k} & k \geq k_0 \ve \omega_{v_k} \neq \phi_{v_k}\\
\phi_{v_k} & k < k_0 \ve \omega_{v_k} \neq \phi_{v_k}
\end{array}
\right.\]
and
\[
\eta_{v_k} = \left\{
\begin{array}{ll}
\omega_{v_k} & \omega_{v_k} = \phi_{v_k} \\
\phi_{v_k} & k \geq k_0 \ve \omega_{v_k} \neq \phi_{v_k}\\
\omega_{v_k} & k < k_0 \ve \omega_{v_k} \neq \phi_{v_k}.
\end{array}
\right.\]

\noindent By the property of our labeling,
\begin{equation} \label{eq:prop_label}
\mu[\sigma]\mu[\eta]\leq
\mu[\sigma^{(k-1)}]\mu[\phi] e^{4 \cE(G) \beta}.
\end{equation}
and $K[\sigma^{(k-1)} \to \sigma^{(k)}] \geq \exp(-2\Delta\beta).$
Now a short calculation concludes the proof:
\begin{eqnarray} \nonumber
\rho &\leq&
\sup_e \sum_{\sigma,\eta \hbox{ s.t. } e \in \gamma(\sigma,\eta)}
 \frac{\mu[\sigma] \mu[\eta]} {\mu[\sigma^{(k-1)}] K[\sigma^{(k-1)} \to
\sigma^{(k)}]} \\
&\leq&
e^{4  \cE(G) \beta}
\sup_e\sum_{\phi} \frac{\mu[\sigma^{(k-1)}] \mu[\phi]}{\mu[\sigma^{(k-1)}]
K[\sigma^{(k-1)} \to \sigma^{(k)}]}
\\ &\leq& \label{eq:canonical}
e^{4  \cE(G) \beta}
e^{2 \Delta \beta} \sum_{\phi} \mu[\phi]
\leq e^{(4  \cE(G) + 2 \Delta) \beta}\,.
\end{eqnarray}
The last inequality follows from the fact
that the map
$\gamma \to \phi$ is injective and therefore $\sum_{\phi} \mu[\phi] \leq 1$.
\end{proof}
\begin{proof}[Proof of part \ref{prop:exp:color} of Proposition \ref{prop:cut-width}]
The previous argument does not directly extend to
 coloring, because the
configurations $\sigma^{(k)}$ along the path (as defined above) may not be
proper colorings.
Assume that $q \geq \Delta + 2$ and let $v_1 < v_2 \cdots < v_{\Vol}$ be an ordering
of the vertices of $G$ which achieves the cut-width. We construct a path
$\gamma(\sigma,\eta)$ such that
\begin{equation} \label{eq:coloring1}
|\gamma(\sigma,\eta)| \leq (\Delta + 1) n.
\end{equation}
Moreover, for all $\tau \in \gamma(\sigma,\eta)$ there exists a $k$ such that
\begin{equation} \label{eq:coloring2}
\tau_v = \left\{ \begin{array}{ll} \eta_v   & \mbox{if } v \leq v_k \\
                                  \sigma_v & \mbox{if } v > v_k \mbox{ and }    v \nsim \{v,\ldots,v_k\}
                \end{array} \right.
\end{equation}
The way to construct a path $\gamma(\sigma,\eta)$ satisfying (\ref{eq:coloring1}) and (\ref{eq:coloring2})
is the following:
$\sigma^0=\sigma$. Given $\sigma^k$, we proceed to create
$\sigma^{k+1}$ as follows: Let
$i(k)=\inf\{j:\sigma^k_{v_j}\neq\eta_{v_j}\}$. If
\[
\rho = \left\{ \begin{array}{ll} \sigma^k_v   & \mbox{if } v \neq v_{i(k)} \\
                                  \eta_v & \mbox{if } v =  v_{i(k)}
                \end{array} \right.
\]
is a legal configuration, then $\sigma^{k+1}=\rho.$ otherwise, let
\[
h(k)=\inf\{j:\sigma^k_{v_j}=\eta_{v_{i(k)}} \ve v_j\sim v_{i(k)}\},
\]
and let $c$ be a color that is different from $\eta_{v_{i(k)}}$ and is
legal for $v_{h(k)}$ under $\sigma^k$. Such a color exists because $q
\geq \Delta + 2$.
Then, we take
\[
\sigma^{k+1}_v = \left\{ \begin{array}{ll} \sigma^k_v   & \mbox{if } v \neq v_{h(k)} \\
                                  c & \mbox{if } v =  v_{h(k)}
                \end{array} \right.
\]

It is easy to verify that the path satisfies (\ref{eq:coloring1}) and (\ref{eq:coloring2}).
Since all legal configurations have the same weight,
(\ref{eq:prop_label}) is replaced by
\begin{equation} \label{eq:prop_label_color}
\mu[\sigma]\mu[\eta]=
\mu[\sigma^{(k-1)}]\mu[\phi]
\end{equation}

On the other hand, the map $\gamma \to \phi$ is not injective. Instead, by (\ref{eq:coloring2}),
there are at most $(q-1)^{\cE(G)}$ paths which are mapped to the same coloring.
We therefore obtain that for the coloring model $\rho \leq n (q-1)^{\cE(G) + 1}$ and therefore from
(\ref{eq:coloring1}) and
(\ref{eq:prop_label_color}),

\[
\tau_2 \leq (\Delta + 1) n q (q-1)^{\cE(G)}.
\]

\end{proof}

\subsection{Hyperbolic graphs} \label{subsec:hyper}
In this subsection we show that balls in a hyperbolic tiling have
logarithmic cut-width. Let $G=(V,E)$ be an infinite planar graph
and let $o\in V$. Let $G_r$ be the ball of radius $r$ in $G$
around $o$, with the induced edges. The following proposition
implies
 Propositions \ref{prop:hyptylexp} and
\ref{prop:corrhyp}.

\begin{Proposition}\label{prop:hyper}
\begin{enumerate}
\item\label{hypecut-width}
Suppose $G$ has 
\begin{itemize}
\item
a positive Cheeger constant, 
\item
degrees bounded by $\Delta$, 
\item
for all $r$, no cycle from $G_r$ separates two
vertices of $G\setminus G_r$,
\end{itemize}
then there
exists constants $\alpha_1$ and $\alpha_2$ s.t. $\cE(G_r)\leq
\alpha_1\Delta \log(|G_r|) + \alpha_2 \Delta.$
\item\label{hypecorrelation}
Assume that $G$ has bounded degrees, bounded co-degrees,
no cycle from $G_r$ separates two vertices of
$G\setminus G_r$, and the following weak
isoperimetric condition holds:
\begin{equation}\label{logcheeger}
|\partial A|\geq C\log(|A|)
\end{equation}
 for every finite
$A\subseteq G$ and for some constant $C$.

Then there exist $\beta'<\infty$ and $\delta>0$ s.t. for every $\beta >
\beta'$, for every $r$ and for
every $u,v$ in $G_r$, the free Gibbs measure for the Ising model
on $G_r$ with inverse temperature $\beta$ satisfies
$\cov(\sigma_u,\sigma_v)\geq\delta.$
\end{enumerate}
\end{Proposition}

\begin{proof}[Proof of part \ref{hypecut-width} of Proposition
  \ref{prop:hyper}]
\newcommand{\maxdeg} {{\Delta}}

Consider a planar embedding of $G$. Since no cycle from $G_r$
separates two vertices of  $G\setminus G_r$, all vertices of $G\setminus G_r$ are
in the same face of $G_r$, and without loss of generality we can
assume that it is the infinite face of our chosen embedding of
$G_r$.

Let $T$ be a shortest path tree from $o$ in $G_r$. In other words, $T$ 
is a tree such that
for every vertex $v$, the path from $o$ to $v$ in $T$ is a shortest
path in $G_r$. Let $e_1\in T$ be an edge adjacent to $o$.
We perform a depth-first-search traversal of $T$, starting from
$o=v_0$, traversing $e_1$ to its end vertex $v_1$, and continuing in
counterclockwise order around $T$. This defines a linear ordering
$v_0\leq v_1\leq \cdots \leq v_{n-1}$ of the vertices of $G_r$.

Consider the induced ordering $w_1\leq w_2\leq\cdots\leq w_k$ on the
vertices of $G_r$ which are at distance exactly $r$ from $o$.

Fix $i<j$. We first consider edges between
\begin{equation*}
V_{ij}=\{ u\in G_r:w_i< u< w_j , u \hbox{ not ancestor of $w_j$ in }T\}
\end{equation*}
and $G_r\setminus V_{ij}$. Note that $V_{ij}$ does not contain any vertex on the
paths in $T$ from $o$ to either $w_i$ or $w_j$. Obviously, there
can be edges from $V_{ij}$ to vertices on the paths from $o$ to $w_i$ or
$w_j$. Let $e=\{ u,v\}$ be an edge with one endpoint in $V_{ij}$ and
the other end in $G\setminus G_r$ but not on Path$(w_i)$ or
Path$(w_j)$,
where
$\hbox{Path}(w_j)$ denote the path in $T$ from $o$ to $w_j$. 
Without loss of generality, assume that $w_i<u<w_j<v$. The case where 
$v < w_i < u < w_j$ is treated similarly.

The path from $o$ to $u$ in $T$, followed by the edge $e$, followed by
the path from $v$ to $o$ in $T$, defines a cycle $C_e$ in $G_r$ (see Figure \ref{fig:hyppath}).
Since $w_i<u<w_j<v$, $C_e$ must enclose exactly one of $w_j$ and
$w_i$. Since the graph is embedded in the plane, the ones among those
cycles which enclose $w_j$ must form  a nested sequence, and
therefore there is an outermost such cycle $C_{e^*}$ with an
associated ``outermost'' edge $e^*=\{ u^*,v^* \}$.
Similarly, among the edges such that the corresponding cycle
encloses $w_i$, there is an ``outermost'' edge $f^*=\{ x^*,y^*\}$.

\begin{figure}[h]
\center{\epsfig{figure=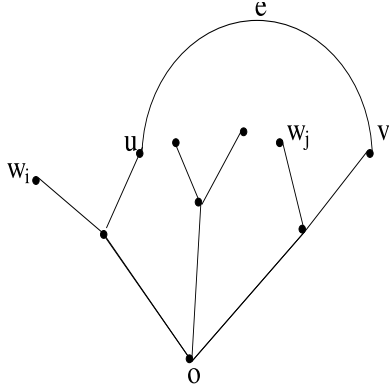, height=2in, width=2in}}
\caption{The cycle $C_e$ defined by $e$.
}\label{fig:hyppath}
\end{figure}


There can only be edges from $V_{ij}$ to the vertices
enclosed by $C_{e^*}$ or by $C_{f^*}$ (note that this includes the
paths from $o$ to $w_i$ and to $w_j$) . Since all the vertices of
$G\setminus G_r$ are in the infinite face of $G_r$, hence outside
$C_{e^*}$, the set of vertices enclosed by $C_{e^*}$ is the same
in $G$ as in $G_r$. Let $A$ denote the set of vertices enclosed by
$C_{e^*}$ (including $C_{e^*}$). We have: $|\partial A|\leq 2r+1$,
hence $|A|\leq (2r+1)/c$, where $c$ is the Cheeger constant of $G$.
Reasoning similarly for $C_{f^*}$, we obtain that the set of
vertices in $G_r\setminus V_{ij}$ adjacent to $V_{ij}$ has size at
most $(4r+2)/c$.

Let $B_j=V_{j-1,j}$ for $j\geq 2$, and $B_1=\{ u\in G_r: u<w_1
\hbox{ or } u>w_k, u\hbox{ not an ancestor of }w_1 \hbox{ in }T\}$. 
Let us bound the cardinality of $B_j$.
As above, we define $C_{e^*}$ and $C_{f^*}$. Let $A$ denote the
union of $B_j$, of the vertices enclosed by $C_{e^*}$, and of the
vertices enclosed by $C_{f^*}$. 

Since the vertices of $B_j$ are at
distance at most $r-1$ from $o$, they have no neighbors in
$G\setminus G_r$. Thus the neighborhood of $A$ in $G$ is such that
$\partial A\subset C_{e^*}\cup C_{f^*}$,
hence $|\partial A|\leq 4r+2$, and so $|B_j|\leq |A|\leq (4r+2)/c$.

Finally, to compute the cut-width, let $S=\{ u: v_0\leq u\leq v_i
\}$, and let $j$ be such that $w_{j-1}<v_i\leq w_j$. We have:
$$V_{1,j-1}\subseteq S \subseteq B_1\cup V_{1,j-1}\cup B_j\cup\hbox{Path}(w_1)\cup
\hbox{Path}(w_{j-1})\cup\hbox{Path}(w_j).$$ Thus
the set of edges between $S$ and $G_r\setminus S$ has size at most
$\Delta (4r+2)/c+ (|B_1\cup B_j|)\Delta + (3r+1)\Delta$, which is at most
$(3r+1+(12r+6)/c)\Delta$. 

Since $G$ has positive Cheeger constant, $$|G_r|=|G_{r-1}\cup\partial
G_{r-1}|\geq |G_{r-1}|(1+c),$$
and so $|G_r|\geq (1+c)^r$, that is, $r\leq \log |G_r| / \log
(c+1)$. Hence the set of edges between $S$ and $G\setminus S$ has size
at most $(3+12/c)(\Delta/\log (c+1))\log |G_r| + (1+6/c)\Delta$.
This concludes the proof.
\end{proof}

\begin{figure}[h]
\center{\epsfig{figure=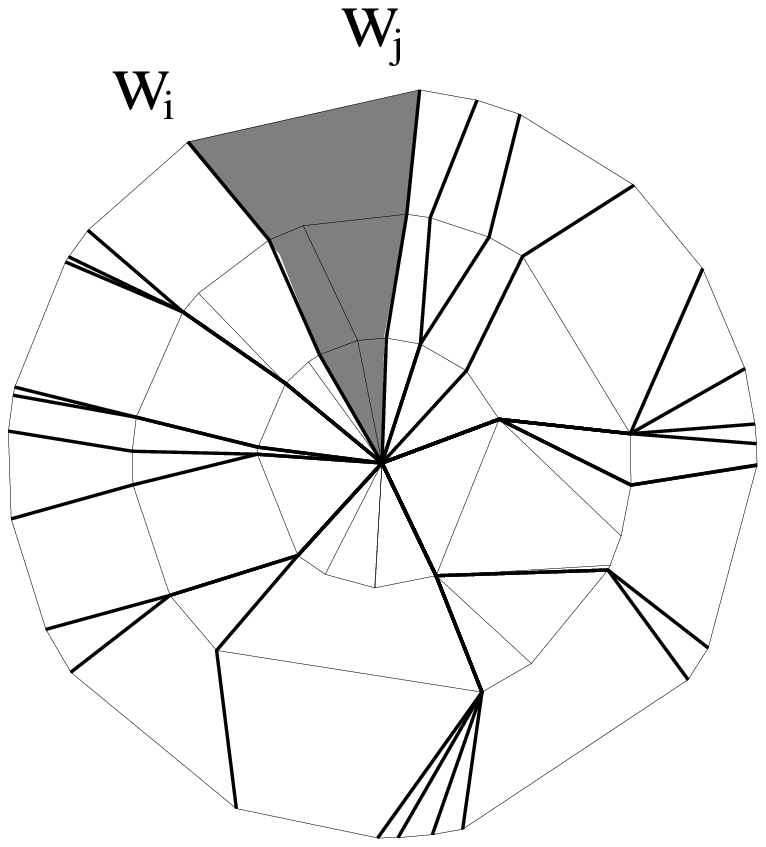, height=1.5in, width=1.5in}}
\caption{The region between $\mbox{Path}(W_i)$ and $\mbox{Path}(W_j)$}\label{fig:hyp}
\end{figure}

\begin{proof}[Proof of part \ref{hypecorrelation} of Proposition \ref{prop:hyper}]
We use the Random Cluster representation of the Ising model (see,
e.g. \cite{FK}) and a standard \prls { }path-counting argument.
For every $u$ and $v$ in $G_r$,
$\cov(\sigma_u,\sigma_v)$ is the probability that $u$ is connected to
$v$ in the Random Cluster model. Fix $p<1$. 
The exact value of $p$ will be specified later.

Then, $\beta$ is large enough, i.e., if 
$(1-e^{-\beta})/e^{-\beta} > 2 p / (1-p)$, 
then the Random Cluster model dominates percolation with parameter
$p$. So, what we need to show is that for a graph satisfying the
requirements of part \ref{hypecorrelation} of the proposition and $p$
high enough, there exists $\delta>0$ s.t. for every $r$ and every $u,v$
in $G_r$, we have $\P_p(u\leftrightarrow v)\geq \delta$. By the FKG inequality
(see \cite{FKG}),
\[
\P_p(u\leftrightarrow v)\geq\P_p(u\leftrightarrow o)\P_p(v\leftrightarrow o)
\]
where $o$ is the center. Therefore we need to show that
$\P(v\leftrightarrow o)$ is bounded away from zero. To this end, we
will pursue a standard path counting technique: in order for $o$ and
$v$ not to be connected, there needs to be a closed path in the dual
graph that separates $o$ and $v$.
\begin{Claim}\label{afikoman}
There exists $M=M(G)$ s.t. for every $r$ and $v\in G_r$ there are at
most $M^k$ paths of length $k$ in the dual graph of $G_r$
 that separate $o$ from $v$.
\end{Claim}
By Claim \ref{afikoman}, if we take $p>1-{1}/{(2M)}$ and choose
$\beta$ accordingly then the probability that there exists
a closed path in the dual graph that separates $o$ and $v$
is bounded away from $1$.
\end{proof}
\begin{proof}[Proof of Claim \ref{afikoman}]
\newcommand {\uu} {{\psi}}
\newcommand{\hatdelta} {{\hat{\Delta}}}
Here again, we consider an embedding of $G_r$ such that
all the vertices of $G\setminus G_r$ lie on the infinite face $F$ of $G_r$.

Let $\gamma$ be a shortest path connecting $v$ to $o$ in $G_r$.
Every dual path separating $v$ from $o$ must intersect $\gamma$.
For an edge $e$ let $\Lambda_k(e)$ be the set of dual paths $\uu$ of
length $k$ separating $o$ from $v$ such that $\uu$ intersects $e$.
If $\hatdelta$ is the maximal co-degree in $G$ then $|\Lambda_k(e)|\leq
\hatdelta^k$
for every $e$.

Let $e\in \gamma$ be such that $d(e,o)>\exp (k/C)+k\hatdelta$,
$d(e,v)> \exp (k/C)+k\hatdelta$, and $d(e,F)>k\hatdelta$. We will now show that
$|\Lambda_k(e)|=0$. Assume, for a contradiction, that $\psi\in
\Lambda_k(e)$. Since $\psi$ has length $k$ and $d(e,F)>\hatdelta k$, $\psi$
does not touch the outer face $F$, and so the area enclosed by $\psi$
in $G_r$ equals the area enclosed by $\psi$ in $G$. The dual path $\psi$ encloses
either $v$ or $o$. Without loss of generality, assume that it encloses $o$.
Let $e'$ be the edge of $\gamma$ closest to $o$ which $\psi$
intersects. Since $\psi$ has length $k$, we get that $d(e',o)>\exp (k/C)$, and so at least
$1+\exp (k/C)$ vertices of $\gamma$ are enclosed by $\psi$. By
(\ref{logcheeger}) this implies that $\psi$ has length strictly greater than $k$,
a contradiction.

Thus, the total number of paths of length $k$ separating $o$ from $v$ is at most
$$\sum_{e:|\Lambda_k(e)|\neq 0} |\Lambda_k(e)|\leq [2\exp(k/C)+k\hatdelta)+
\hatdelta k]\hatdelta^k.$$


\end{proof}

\begin{figure}[h]
\center{\epsfig{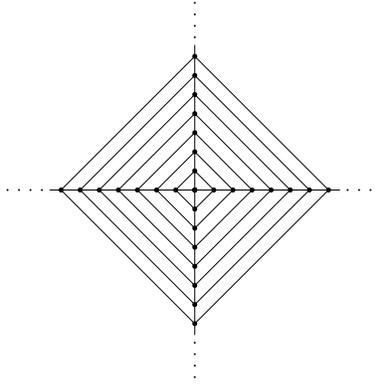}}
\caption{An example of a planar graph s.t. for all temperatures, correlations decay
exponentially with distance.
}\label{fig:diamond}
\end{figure}

{\bf Remark:} An isoperimetric inequality of the type of (\ref{logcheeger}) is necessary. An example where all other conditions of part \ref{hypecorrelation} of Proposition \ref{prop:hyper} are satisfied and yet the conclusion does not hold can be found in Figure \ref{fig:diamond}.

\subsection{A polynomial upper bound for trees} \label{subsec:recur}
In this subsection
we give an improved bound on relaxation time for the tree.

Let $A$ be a finite set, and let $\alpha_{vw}:A\times A\to\R_{+}$ be
a weight function. Let $G$ be a graph. Let the Glauber dynamics be as
defined above, and let $\L=\L(A, \alpha, G)$ be its generator. We say
that the Glauber dynamics on $(A, \alpha, G)$ is ergodic if for every
two legal configurations $\sigma_1$ and $\sigma_2$, we have
$\left(\exp(\L)\right)_{\sigma_1\sigma_2}>0$.
We will prove the following proposition:
\begin{Proposition}\label{prop:polbound}
Let $b\geq 2$, and let $T$ denote the infinite $b$-ary tree, and let $T_n$
be the $b$-ary tree with $n$ levels.  If the Glauber dynamics on
$(A, \alpha, T_n)$ is ergodic for every $n$ then
\begin{equation*}
\limsup_{n\to\infty}\frac{1}{n}\log\left(\tau_2\left(\L\left(A,
\alpha, T_n\right)\right)\right) <\infty .
\end{equation*}
\end{Proposition}

\begin{Conjecture}\label{conj:pol}
Let $b\geq 2$, $T$ denote the infinite $b$-ary tree, and let $T_n$
be the $b$-ary tree with $n$ levels. If the Glauber dynamics on
$(A, \alpha, T)$ is ergodic then there exists $0\leq\tau<\infty$
s.t.
\begin{equation}\label{eq:conjpoll}
\lim_{n\to\infty}\frac{1}{n}\log\left(\tau_2\left(\L\left(A,
\alpha, T_n\right)\right)\right) = \tau .
\end{equation}
\end{Conjecture}

We prove a special case of Conjecture \ref{conj:pol}:
\begin{Proposition}\label{prop:soft}
If the interactions are soft, i.e. $\alpha_{vw}(a,b) > 0$ for all
$v,w,a$ and $b$,
then (\ref{eq:conjpoll}) holds.
\end{Proposition}

The main tool we use for proving Propositions \ref{prop:polbound}
and \ref{prop:soft} is block dynamics (see e.g. \cite{Ma}). For a
spin (or a color) $a\in A$, we denote by $\L(a,\alpha,n)$ the
Glauber dynamics on the $b$-ary tree of depth $n$, under the
interaction matrix $\alpha$ and with the boundary condition that
the root has  a parent colored $a$. With a slight abuse of
notations, we say that $\tau_2(a,\alpha,n)$ is the relaxation
time for $\L(a,\alpha,n)$.

\newcommand {\htt} {{\hat{\tau}}}
\begin{Lemma}\label{subadd}
Let
\begin{equation*}
\htt_2(\alpha,n)=\sup_{a\in A}\tau_2(a,\alpha,n)
\end{equation*}
Then,  for all $m$ and $n$,
\begin{equation*}
\htt_2(\alpha,n+m)\leq\htt_2(\alpha,n)\htt_2(\alpha,m).
\end{equation*}
\end{Lemma}

\begin{proof}
Let $l=n+m$. Partition the tree $T_l$ into disjoint sets
$V_1,...,V_k$ to be specified below. We call  $V_1,...,V_k$ {\em
blocks}, and consider the following block dynamics: Each block
$V_i$ has a (rate $1$) Poisson clock, and whenever it rings, $V_i$
updates according to its Gibbs measure determined by the boundary
conditions given by the configurations of $T^{(b)}_l-V_i$ and by
the external boundary conditions. We denote by
$\L^B=\L^B(V_1,...,V_k)$ the generator for the block dynamics, and
let $\L^B_a$ be the generator for the block dynamics with the
boundary condition that the parent of the root has color $a$.

By \cite{Ma}[Proposition 3.4, page 119],
\begin{equation*}
\htt_2(\alpha,l)\leq\sup_{i}\htt_2(\alpha,V_i)\cdot\sup_{a\in
  A}\tau_2(\L^B_a)
\end{equation*}
We now define the partition to blocks. For every vertex $v$ up to
depth $n$, the singleton $\{v\}$ is a block, and for every vertex
$w$ at depth $n$, the full subtree of depth $m$ starting at $w$ is
a block (see Figure \ref{blocks}).
\begin{figure}[h]
\center{\epsfig{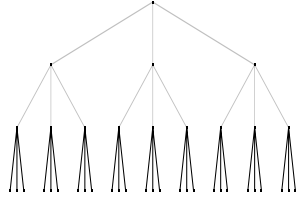}}
\caption{Partition of a tree to blocks}\label{blocks}
\end{figure}
All we need now to finish the proof is the following easy claim:
\begin{Claim}\label{blocdin}
\begin{equation*}
\sup_{a\in  A}\tau_2(\L^B_a) = \htt_2(\alpha,n).
\end{equation*}
\end{Claim}
\begin{proof}
We use the following fact (that could also serve as a definition of
the relaxation time). Given the dynamics $\L$ we define the Dirichlet
form $\dir[g,g] =\frac{1}{2}
\sum_{\sigma,\tau} \mu[\sigma] K[\sigma \to \tau] (g(\sigma) -
g(\tau))^2$. Then
\begin{equation}\label{eq:dir_def}
\tau_2=
\sup\left\{ \frac
{\mu[g^2]}
{\dir[g,g]}
~:~ \mu[g]=0 \right\}.
\end{equation}
Clearly, the expression in (\ref{eq:dir_def}) evaluated for $f$
and $\L^B_a$ is equal to the one evaluated for $g$ and
$\L(a,\alpha,n)$, if
\begin{equation}\label{elchanan}
g(\eta)=f(\sigma) \mbox{ for all } \eta \mbox{  and } \sigma \mbox{
  s.t. } \eta\left| _{T_n}\right.=\sigma.
\end{equation}
Therefore, we need to show that the maximum in (\ref{eq:dir_def})
for the dynamics $\L^B_a$ is
obtained at a function that satisfies (\ref{elchanan}). The maximum in
(\ref{eq:dir_def}) is obtained at an eigenfunction of
$\L^B_a$. Moreover for
every function $g$, $\L^B_a(g)$ satisfies (\ref{elchanan}) with some
function $f$. It now follows that the maximum is obtained at a
function that satisfies (\ref{elchanan}).
\end{proof}
\end{proof}

\begin{proof}[Proof of Proposition \ref{prop:polbound}]
From Lemma \ref{subadd} and the sub-additivity lemma, we learn that
\begin{equation*}
\limsup_{n\to\infty}\frac{1}{n}\log\left(\htt_2\left(\L\left(A, \alpha, T_n\right)\right)\right)
<\infty
\end{equation*}
By another application of Matinelli's block dynamics lemma, we get
that
\begin{equation}\label{gadolkatan}
\tau_2\left(\L\left(A, \alpha, T_n\right)\right)\leq
\tau_2\left(\L\left(A, \alpha, T_1\right)\right) \cdot
\htt_2\left(\L\left(A, \alpha, T_{n-1}\right)\right)
\end{equation}
and the proposition follows.
\end{proof}
\begin{proof}[Proof of proposition \ref{prop:soft}]
From Lemma \ref{subadd} and the sub-additivity lemma, we learn that
there exists $0\leq\tau<\infty$ s.t.
\begin{equation*}
\lim_{n\to\infty}\frac{1}{n}\log\left(\htt_2\left(\L\left(A, \alpha,
      T_n\right)\right)\right) = \tau.
\end{equation*}
For every $a$, let $\mu_a$ be the Gibbs measure for the tree of
depth $n$ with the boundary condition that the parent of the root
has color $a$. Note that $\mu_a$ is the stationary distribution of
$\L(a,\alpha,n)$. Since the interactions are soft, there exists
$0<C<\infty$ s.t. for every $a$,  every $n$, and every two
configurations $\sigma$ and $\eta$ on the tree of depth $n$,
\[
\frac{1}{C}\mu(\sigma) \leq \mu_a(\sigma) \leq C\mu(\sigma),
\]
and
\[
\frac{1}{C}\L(\alpha,n)_{\sigma,\eta} \leq
\L(a,\alpha,n)_{\sigma,\eta} \leq C\L(\alpha,n)_{\sigma,\eta}.
\]
Therefore, by (\ref{eq:dir_def}),
\[
\frac{1}{C^3}\tau_2\left(\L\left(A, \alpha,T_n\right)\right) \leq
\htt_2\left(\L\left(A, \alpha,T_n\right)\right) \leq
C^3\tau_2\left(\L\left(A, \alpha,T_n\right)\right)
\]
and the proposition follows.
\end{proof}

\section{Lower Bounds}\label{sec:lowerbounds}
\begin{proof}[Proof of Theorem~\ref{thm:tree}, part~\ref{treesuperlinear}]
Theorem~\ref{thm:tree}
part~\ref{treesuperlinear} is a direct
consequence of the extremal characterization of $\tau_2$ given
in (\ref{eq:dir_def}),
applied to the particular test function $g$
which sums the spins on the boundary
of the tree.
It is easy to see that $\mu[g] = 0$ and that
\[
\dir[g,g] \leq \sum_{\sigma,\tau} \mu[\sigma] K[\sigma \to \tau] = O(n_r).
\]
We repeat the variance calculation from ~\cite{EKPS}. When $b(1 - 2
\eps)^2 > 1$:
\begin{eqnarray*}
\mu[g^2]&=&
\sum_{w\in\partial\T}\mu[\sigma_w^2]+
\sum_{w\in\partial\T}\sum_{{v\in\partial\T \atop v\neq w}}
\mu[\sigma_w\sigma_v]\\
&=&b^r\cdot\left(
1+ \sum_{i=1}^r (b-1) b^{i-1} (1-2\eps)^{2i} \right)\\
&=& b^r \left(1+\Theta\left(\left(b(1-2\eps)^2\right)^r\right)\right)\\
&=& \Theta\left(n_r^{1+\log_b(b(1-2\eps)^2)}\right).
\end{eqnarray*}
It now follows by (\ref{eq:dir_def})
that if $b(1-2\eps)^2 > 1$ then
\[
\tau_2=\Omega \left(n_r^{\log_b(b(1-2\eps)^2)}\right),
\]
as needed.
Repeating the calculation for the case $b (1 - 2 \eps)^2 = 1$ yields that 
\[
\tau_2 = \Omega(\log n_r).
\]
The proof follows.
\end{proof}

\begin{Remark}
Suppose that $\mu$ admits a Markovian representation where the
conditional distribution of $\sigma_u$ given its parent $\sigma_v$ is
given by an $|\A| \times |\A|$ mutation matrix $P$. Let $\lam_2(P)$ be
the second eigen-value of $P$ (in absolute value), and $x$ the
corresponding eigen-vector, so that $Px^t = \lam_2(P) x^t$ and $|x|_2 = 1$.

Let $g$ be the test function $g = c_n x^t$, where $c_n(i)$ is the
number of boundary nodes that are labeled by $i$.
It is then easy to see once again that $\dir[g,g] = O(n_r)$.
Repeating the calculation from \cite{MP} it follows that if $b
|\lam_2(P)|^2 > 1$, then
\[
\Var[g] = \Theta\left(n_r^{1+\log_b(b |\lam_2(P)|^2)}\right).
\]
Thus in this case,
\[
\tau_2 = \Omega \left(n_r^{\log_b(b |\lam_2(P)|^2)}\right).
\]
\end{Remark}

\begin{figure} [hptp] \label{rec_maj}
\epsfxsize=5cm
\hfill{\epsfbox{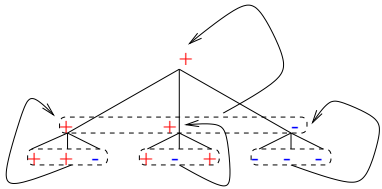}\hfill}
\caption{The recursive majority function.}
\end{figure}

In order to prove the lower bound on the relaxation time for
very low temperatures stated in Theorem~\ref{thm:tree}
part~\ref{treeinfinitedegree},
we apply (\ref{eq:dir_def}) to the test function $g$ which is obtained
by applying recursive majority to the boundary spins;
see \cite{Mo2} for background regarding the recursive-majority function for
the Ising model on the tree.
For simplicity we consider first the ternary tree $T$, see
Figure $5$.
Recursive majority is
defined on the configuration space as follows.
Given a configuration $\sigma$, first label each boundary vertex $v$
by its spin $\sigma_v$.
Next, inductively label each interior vertex $w$ with the label of the
majority of the children of $w$.
The value of the recursive majority function $g$ is then the label of
the root. We write $\sigma_v$ for the spin at $v$ and $m_v$ for the
recursive majority value at $v$.

\begin{Lemma}
If $u$ and $w$ are children of the same parent $v$, then
$\P[m_u \neq m_w] \leq 2 \eps + 8 \eps^2$.
\end{Lemma}
\noindent {\bf Proof:}
$$
\P[m_u \neq m_w] \leq  \P[\sigma_u \neq m_u] + \P[\sigma_w \neq
m_w] + \P[\sigma_u \neq \sigma_v] +
            \P[\sigma_w \neq \sigma_v].
$$

 \noindent We will show that recursive majority is highly correlated with spin, {\em
i.e.}
if $\eps$ is small enough (say $\eps < 0.01$),
then $\P[m_v \neq \sigma_v] \leq 4 \eps^2$.

The proof is by induction on the distance $\ell$ from $v$ to the boundary of
the tree.
For a vertex $v$ at distance $\ell$ from the boundary of the tree, write
$p_\ell = \P[m_v \neq \sigma_v]$.
By definition $p_0 = 0 \leq 4 \eps^2$.

For the induction step, note that if $\sigma_v \neq m_v$ then one of the
following events hold:
\begin{itemize}
\item
At least $2$ of the children of $v$, have different $\sigma$ value than that
of $\sigma_v$, or
\item
One of the children of $v$ has a spin different from the spin at $v$, and
for some other child
$w$ we have $m_w \neq \sigma_w$, or
\item
For at least $2$ of the children of $v$, we have $\sigma_w \neq m_w$.
\end{itemize}
Summing up the probabilities of these events, we see that $
p_\ell \leq 3 \eps^2 + 6 \eps p_{\ell-1} + 3 p_{\ell-1}^2.$
It follows  that $p_\ell \leq 4 \eps^2$, hence the Lemma.
$\qed$

\begin{proof}[Proof of Theorem~\ref{thm:tree} part~\ref{treeinfinitedegree}]
Let $m$ be the recursive majority function. Then by symmetry $\E[m] = 0$,
and $\E[m^2] = 1$.
By plugging $m$ in definition (\ref{eq:dir_def}), we see that
\begin{equation} \label{eq:rec_maj_cut}
\tau_2 \geq \left( \sum_{\sigma, \tau : m[\sigma] = 1, m[\tau] = -1}
\mu[\sigma] \P[\sigma \to \tau] \right)^{-1}.
\end{equation}
Observe that if $\sigma,\tau$ are adjacent configurations
(i.e.,  $\P[\sigma \to \tau] > 0$) such that
$m(\sigma) = 1$ and $ m(\tau) = -1$,
then there is a unique vertex $v_r$ on the boundary of the tree
where $\sigma$ and $\tau$ differ.
Moreover, if $\rho=v_1,\ldots,v_r$ is the path from $\rho$ to $v_r$, then
for $\sigma$ we have
$m(v_1) = \ldots = m(v_r) = 1$ while for $\tau$ we have $m(v_1) = \ldots =
m(v_r) = -1$.
Writing $u_i,w_i$ for the two siblings of $v_i$ for $2 \leq i \leq k$, we
see that for
all $i$, for both  $\sigma$ and $\tau$ we have $m(u_i) \neq m(v_i)$.
Note that these events are independent for different values of $i$.
We therefore obtain that the probability that $v_1,\ldots,v_r$ is such a
path is bounded by
$(2 \eps + 8 \eps^2)^{r-1}$. Since there are $3^r$ such paths and
since $\P[\sigma \to \tau] \leq 3^{-r}$ we obtain
that the right term of (\ref{eq:rec_maj_cut}) is bounded below by
\[
(2 \eps + 8 \eps^2)^{1-r} \ge \sT^{\Omega( \beta )}\, .
\]
\end{proof}
Note that the proof above easily extends to the $d$-regular tree
for $d \geq 3$. A similar proof also applies to the binary tree
$T$, where $g$ is now defined as follows. Look at $T_k$ for even
$k$. For the boundary vertices define $m_v = \sigma_v$. For each
vertex $v$ at distance $2$ from the boundary, choose three leaves
on the boundary below it $v_1,v_2,v_3$ and let $m_v$ be the
majority of the values $m_{v_i}$. Now continue recursively.

Repeating the above proof, and letting
$p_{\ell} = P[m_v \neq \sigma_v]$ for a vertex at distance $2 \ell$,
we derive the following recursion:
$p_{\ell} \leq 3 (2 \eps)^2 + 6 (2 \eps) p_{\ell - 1} + 3 p_{\ell - 1}^2$.
We then continue in exactly the same way as for the ternary tree.

\section{High temperatures}\label{sec:hightemp}

\begin{proof}[Proof of Theorem~\ref{thm:tree} part~\ref{treelinear}]
Our analysis uses a comparison to block dynamics.

{\bf Block dynamics}.
We view our tree $T = T_r^{(b)}$ as a part of a larger $b$-ary tree $T_*$ of
height
$r+2h$, where the root  $\rho$ of $T$ is at level $h$ in $T_*$.
For each vertex $v$ of $T_*$, consider
the subtree of height $h$ rooted at $v$.
A {\bf block} is by definition the intersection of
$T$ with
such a subtree.
Each block has a rate $1$ Poisson clock and whenever the clock rings
we
erase all the spins
of vertices belonging to the block, and put new spins in, according to the
Gibbs distribution conditional on the spins in the rest of $T$.

\noindent
{\bf Discrete dynamics:}
In order to be consistent with \cite{BD}, we will first analyze the
corresponding discrete time dynamics: at each step of the block
dynamics, pick a block at random, erase all the spins
of vertices belonging to the block, and put new spins in, according to the
Gibbs distribution conditional on the spins in the rest of $T$.

{\bf A coupling analysis}.
We use a weighted Hamming metric on configurations,
\[
d(\sigma,\eta)=
\sum_{v} \lambda^{|v|} 1(\sigma_v\neq\eta_v),
\]
where $|v|$ denotes the distance from vertex $v$ to the root.
Let $\theta=1-2\epsilon$ and $\lambda=1/\sqrt{b}$. Note that
$b \lambda \theta
<1$ and
$\theta < \lambda$.
Starting from two distinct configurations $\sigma$ and $\eta$,
our coupling always  picks the same block in $\sigma$ and in $\eta$ and chooses the
coupling between the two block moves which minimizes $d(\sigma',\eta')$.

We use path-coupling~\cite{BD}, {\em i.e.}, we will prove that for every
pair of configurations
which differ by a single spin, applying one step of the block dynamics will
reduce
the expected distance between the two configurations.

Let $v$ be the single vertex,  such that $\sigma_v\neq \eta_v$.
Then  $d(\sigma,\eta)=\lambda^{|v|}$. Let $B$ denote the
chosen block, and
$\sigma',\eta'$ be the configurations after the move.
In order to understand $(\sigma',\eta')$, we will need the following
Lemma.

\begin{Lemma}\label{reductiontofree}
Let $T$ be a finite tree and let $v\neq w$ be vertices in $T$. Let
$\{\beta_e\geq 0\}_{e\in E(T)}$ be the (ferromagnetic)
interactions on $T$, and let $\{-\infty<H(u)<\infty\}_{u\in V(T)}$
be an external field on the vertices of $T$.
we consider the following conditional Gibbs measures:\\
$\mu_{+,H}$: The Gibbs measure with external field $H$ conditioned on $\sigma_v = 1$. \\
$\mu_{-,H}$: The Gibbs measure with external field $H$ conditioned on $\sigma_v = -1$. \\
Then, the function $\mu_{+,H}[\sigma_w]-\mu_{+,H}[\sigma_w]$ achieves
its maximum at $H\equiv 0$.
\end{Lemma}

Before proving the Lemma, we utilize it 
to prove Theorem \ref{thm:tree}, part \ref{treelinear}.
There are four situations to consider.

\noindent{\bf Case 1.}
if $B$ contains neither $v$ nor any vertex adjacent to $v$, then
 $d(\sigma',\eta')=d(\sigma,\eta)$.

\noindent{\bf Case 2.}
If $B$ contains $v$, then $\sigma '=\eta '$ and
$d(\sigma',\eta')=0=d(\sigma,\eta)-\lambda^{|v|}$.
There are $h$ such blocks, corresponding to the $h$ ancestors
of $v$ at $1,2,\ldots ,h$ generations above $v$. (Note that this holds even
when $v$ is the root of $T$ or a leaf of $T$, because of our definition of
blocks).

\noindent{\bf Case 3.}
If $B$ is rooted at one of $v$'s children, then the conditional
probabilities given the outer boundaries of $B$ are not the same since one block has $+1$ above it
and the other block has $-1$ above it. However both blocks have their
leaves adjacent to the same boundary configuration. When considering the process on the block,
the influence of the boundary configuration can be counted as altering the external field.
Since $\sigma$ and $\eta$ have the same external fields and the same boundary configuration on all
of the boundary vertices except $v$,
by Lemma  \ref{reductiontofree},
conditioning on this lower boundary can only reduce $d(\sigma',\eta')$. Therefore, we bound
$d(\sigma',\eta')$ by studying the case where one block is conditioned to
having
a $+1$ adjacent to the root, the other block is conditioned to having a $-1$
adjacent to the root, and no external field or boundary conditions. Then the
block is
simply filled in a top-down manner, every edge is faithful ({\em i.e.} the
spin of the
current vertex equals the spin of its parent) with probability
$\theta$ and cuts information (the spin of the current vertex is a new
random spin)
 with probability $1 - \theta$. Coupling these
choices for corresponding edges for $\sigma$ and for $\eta$, we see that the
distance between $\sigma'$ and $\eta'$ will be equal to the weight of the
cluster
containing $v$, in expectation
$\sum_j \lambda^{|v|+j} b^j \theta^j
\leq \lambda^{|v|}/(1 - b \lambda \theta)$. There are $b$ such blocks,
corresponding to the $b$ children of $v$.

\noindent{\bf Case 4.}
If $B$ is rooted at $v$'s ancestor exactly $h+1$ generations above
$v$, then the conditional probabilities are not the same since one block has
a leaf $v$ adjacent
to a $+1$ and the other block has a leaf adjacent to a $-1$. There is
exactly one such block.
Again we appeal to Lemma~\ref{reductiontofree} to show that the expected
distance is dominated by
the size of the $\theta$ cluster of $w$.
The expected weight of $v$'s cluster is bounded by summing over the ancestors
$w$ of $v$:
\begin{eqnarray*}
\lefteqn{\sum_{w} \theta^{|v|-|w|} \sum_j \lambda^{|w|+j} b^j \theta^j=}\\
 &=&
\frac{ \sum_{w} \lambda^{|w|} \theta^{|v| - |w|} }{1 - b \lambda \theta}\\
 &=&
 \frac{\lambda^{|v|}}{(1 - \theta \lambda^{-1}) (1 - b \lambda \theta)}.
\end{eqnarray*}

Overall, the expected change in distance is
\[
\begin{aligned}
&\E(d(\sigma',\eta')-d(\sigma,\eta)) \leq \\[1ex]
&\left(\frac{b\lam^{|v|}}{1-b\lambda \theta}
+ \frac{\lambda^{|v|}}{(1 - \theta \lambda^{-1}) (1 - b \lambda
  \theta)} -h\lam^{|v|} \right) \frac{1}{\sT + h - 1}.
\end{aligned}
\]
If the block height $h$ is a sufficiently large constant,
we get that for some positive constant $c$,
\begin{equation}\label{eq**}
\E(d(\sigma',\eta')-d(\sigma,\eta))\leq \frac{- c \lambda^{|v|}}{\sT}
 \leq \frac{-c }{\sT}d(\sigma,\eta).
\end{equation}

Note that $\max d(\sigma,\eta) = \sum_{j \leq r} b^j \lam^j \leq
\sqrt{\sT}$.
Therefore, by a path-coupling argument (see \cite{BD})
we obtain a mixing time of at most $O(\sT \log \sT)$ for the blocks
dynamics.

{\bf Spectral gap of discrete time block dynamics.}
The $(1-c/\sT)$ contraction at each step of the coupling implies, by
an argument from~\cite{Chen} which
we now recall,
that the spectral gap of the block dynamics is at least $c/\sT$.
Indeed, let $\lambda_2$ be the second largest eigenvalue in absolute
value,
and $f$ an eigenvector for $\lambda_2$.
Let $M=\sup_{\sigma,\eta} |f(\sigma)-f(\eta)|/d(\sigma,\eta)$ and
denote by
 $\P$  the transition operator. Then
\begin{eqnarray*}
|\lambda_2|M &=& \!\!\!\!
    \sup_{\sigma,\eta}
    \frac{|\P f(\sigma)-\P f(\eta)|}{d(\sigma,\eta)}~~
    \hbox{ since $f$ eigenvector for $\lambda_2$} \\
&\leq & \!\!\!\!
    \sup_{\sigma,\eta} \sum_{\sigma',\eta'}
    \P[(\sigma,\eta)\rightarrow
    (\sigma',\eta')]\frac{|f(\sigma')-f(\eta')|}{d(\sigma',\eta' )}
    \frac{d(\sigma',\eta')}{d(\sigma,\eta)}\\
&\leq & \!\!\!\!
    \sup_{\sigma,\eta}\sum_{\sigma',\eta'}\P[(\sigma,\eta)\rightarrow
    (\sigma',\eta')]M\frac{d(\sigma',\eta')}{d(\sigma,\eta)}\\
&=& \!\!\!\!
    M \sup_{\sigma,\eta} \frac{\E[d(\sigma',\eta')]}{d(\sigma,\eta)}\\
&\leq & \!\!\!\!
    (1-c/\sT) M ~~\hbox{by ~(\ref{eq**}).}
\end{eqnarray*}
Thus $|\lambda_2| M\leq (1-c/\sT)M$, whence the
(discrete time) block dynamics has relaxation time at most $O(\sT)$.

{\bf Relaxation time for continuous time block dynamics}.  The
continuous time dynamics is $n$ times faster than the discrete time
dynamics. This is true because the transition matrix for the discrete
dynamics is $M=\I+{\frac{1}{n}}\L$ where $\I$ is the $2^n$-dimensional
unit matrix.  Therefore
\[
\tau_2(\mbox{block dynamics})=O(1).
\]

{\bf Relaxation time for single-site dynamics}.
Since each block update can
be simulated by doing a constant number of single-site updates inside the
block,
and each tree vertex only belongs to a bounded number of blocks,
it follows from proposition 3.4 of \cite{Ma} that the  relaxation time
of the single-site Glauber dynamics is also $O(1)$.
\end{proof}

\begin{proof}[Proof of Lemma~\ref{reductiontofree}]

\noindent{\bf Reduction from trees to paths.}
We first claim that it suffices to prove the lemma when the
tree $T$ consists of a path
$v = v_1, \ldots, v_k = w.$
(see Figure~\ref{tree2path}).
To see this, let $T_1,T_2,\ldots,T_k$ be the connected components of $T$
when the edges in the path $v_1,v_2,\ldots,v_k$ are erased,
s.t. $v_i\in T_i$ for $i=1,2,\ldots,k.$
Let $\sigma$ be a configuration on $v_1,\ldots,v_k$, and for a subgraph
$J$ let $S(J)$ be the space of configurations on $J$.
The probability of a configuration $\sigma$ on $v_1,\ldots,v_k$ is
\begin{eqnarray*}
\frac{1}{Z}
&\exp&\left(
\sum_{i=1}^{k-1}\beta_{\{v_i,v_{i+1}\}}\sigma_{v_i}\sigma_{v_{i+1}}
\right)\cdot
\prod_{i=1}^{k}\left(
\sum_{\tau\in S(T_i-\{v_i\})}\exp\left(\hamil(\tau\cup\sigma_{v_i})\right)
\right)\\
=\frac{1}{Z^\prime}
&\exp&\left(
\sum_{i=1}^{k-1}\beta_{\{v_i,v_{i+1}\}}\sigma_{v_i}\sigma_{v_{i+1}}+
\sum_{i=1}^{k}H^\prime_{v_i}\sigma_{v_i}
\right)
\end{eqnarray*}
for some external field $\{H^\prime_u\}$ depending only on $\{H_u\}$ and $\{\beta_e\}$,
where $Z$ and $Z^\prime$ are partition functions and $\hamil(\cdot)$ denotes the Hamiltonian.
\begin{figure}[ht]
\begin{center}{\small
\setlength{\unitlength}{1cm}
\begin{picture}(7,3)
\multiput(0,0)(4.5,0){2}{
  \put(0,0.5){\circle*{.15}}
  \put(1.5,2){\circle*{.15}}
  \put(2,2.5){\circle*{.15}}
  \put(0,0.5){\line(1,1){.5}}
  \put(1.5,2){\line(-1,-1){.5}}
  \put(2,2.5){\line(-1,-1){.5}}
  \put(-0.4,.5){$v_k$}
  \put(1.1,2){$v_2$}
  \put(1.6,2.5){$v_1$}
}

  \put(0,0.5){\line(1,1){.5}}
  \put(1.5,2){\line(-1,-1){.5}}
  \put(2,2.5){\line(-1,-1){.5}}
  \put(2,2.5){\line(1,-1){.5}}
  \put(1.5,2){\line(1,-1){.5}}
  \put(0,0.5){\line(1,-1){.5}}
\put(0,0.5){\line(0,-1){.5}}
\put(1.5,2){\line(0,-1){.5}}
\put(2,2.5){\line(0,-1){.5}}
\put(0,0){\line(1,0){.5}}
\put(1.5,1.5){\line(1,0){.5}}
\put(2,2){\line(1,0){.5}}
\put(0.6,0){$T_k$}
\put(2.1,1.5){$T_2$}
\put(2.6,2){$T_1$}

\put(3,1){\vector(1,0){1}}

\end{picture}
}\end{center}
\caption{Reduction from trees to paths.}\label{tree2path}
\end{figure}
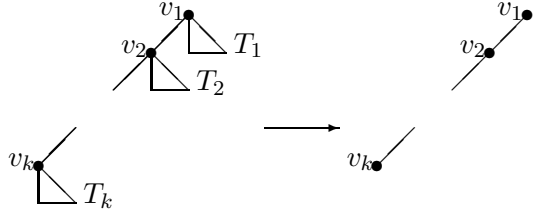

We will now prove the lemma by induction on the length of the path
$v_1,\ldots,v_k$.

\noindent{\bf Paths of length $2$.}
Assume $k=2$. Writing $\beta$ for the strength of $(v_1,v_2)$
interaction,
$H$ for external field at $w = v_2$. Then,
\begin{eqnarray*}
{\mu_{+,\tau}[\sigma_w] - \mu_{-,\tau}[\sigma_w]}
 &=&
   \frac{e^{\beta + H } - e^{-\beta - H }}{e^{\beta + H } +
e^{-\beta - H }} -
   \frac{e^{-\beta + H } - e^{\beta - H }}{e^{-\beta + H } +
e^{\beta - H }}\\
   &=& \tanh (\beta + H ) - \tanh(H  - \beta).
\end{eqnarray*}
It therefore suffices to prove that for $\beta > 0$, the function
\[
H  \mapsto g(\beta,H ) = \tanh (H  + \beta) - \tanh(H  -
\beta)
\]
has a unique maximum at $H =0$.
Consider the partial derivative,
\begin{equation} \label{eq:derg}
g_{H }(\beta,H ) = \cosh^{-2}(H  + \beta) -
                           \cosh^{-2}(H  - \beta).
\end{equation}
Therefore, if $\beta > 0$ and $H  > 0$ then
$g_{H }(\beta,H ) < 0$ and if $\beta > 0$ and $H  < 0$
then
$g_{H }(\beta,H ) > 0$.
Thus $H =0$ is the unique maximum and the claim for $k=2$ follows.

{\bf Induction step.}
We assume that the claim is true for
$k-1$ and prove
it for $k$. We denote $v' = v_{k-1}, \mu'_{+,H} = \mu_{H}[ \cdot | \sigma_{v'} =
1]$ and similarly
$\mu'_{-,H}.$
Now,
\begin{eqnarray}
\lefteqn{\mu_{+,H}[\sigma_w] - \mu_{-,H}[\sigma_w]} \nonumber \\
 &=& \nonumber
\left(
\mu_{+,H}[\sigma_{v'} = 1] \mu'_{+,H}[\sigma_w] + \mu_{+,H}[\sigma_{v'} = -1]
\mu'_{-,H}[\sigma_w]
\right) \\ \nonumber &-&
\left(
\mu_{-,H}[\sigma_{v'} = 1] \mu'_{+,H}[\sigma_w] + \mu_{-,H}[\sigma_{v'} = -1]
\mu'_{-,H}[\sigma_w]
\right)
\\ &=& \label{eq:ind_free}
\frac{1}{2} (\mu_{+,H} - \mu_{-,H})[\sigma_{v'}] (\mu'_{+,H} -
\mu'_{-,H})[\sigma_w].
\end{eqnarray}
Since by the induction hypothesis both multipliers in
(\ref{eq:ind_free}) achieve their maximums at $H\equiv 0$, we get that
${\mu_{+,H}[\sigma_w] - \mu_{-,H}[\sigma_w]}$ also achieves its
maximum at $H\equiv 0$.
\end{proof}

\section{Proof of Theorem ~\ref{thm:general}}\label{sec:general}

Recall that we denoted by $\sigma_r$ the configuration on all vertices
at distance exactly
$r$ from $\rho$. Also recall that $\mu$ is the Gibbs measure
which is stationary for the Glauber dynamics. We abbreviate
$\int f \, d\mu$ as $\mu(f)$.

{\bf Mutual information and $L^2$ estimates.} For Markov chains
such as $\{\sigma_r\}$, it is generally known \cite{Sa} that
(\ref{eq:mut_inf_decay}) follows from  (\ref{eq:cor_decay}), which
itself, is a consequence of the following stronger statement:

There exists $c_*>0$ such that for  any vertex set
$A \subset G_{r/2}$ and  any functions
 $f,g$ with $\mu(f)=\mu(g)=0$, we have
\begin{equation} \label{eq:l2_decay}
\mu(fg) \le  e^{-c_* r} (\mu(f^2)\mu(g^2))^{1/2} \, ,
\end{equation}
provided that $f(\sigma)$ depends only on $\sigma_A$ and $g(\sigma)$ depends
only on $\sigma_r$.
(\ref{eq:l2_decay}) will follow from a more general proposition
below. For a set $A$ of vertices in a graph $G$ we write $\partial_i
A$ for the set of vertices $v$ in $A$ for which there exists an edge
$(v,u)$ with $u \notin A$.

\begin{Proposition}\label{prop:gendecay}
Let $G$ be a finite graph, and let $A$ and $B$ be sets of vertices in $G$.
Let $d$ be the distance between $A$ and $B$ and let $\Delta$ be the
maximum degree in $G$. For $0 < c < 1$, let
\begin{equation}\label{eq:largedivrate}
I(c)=c-\log c-1.
\end{equation}
Let $c^{\ast}$ be the unique $0 < c < 1 $ satisfying $I(c) = \log
\Delta$ and for $0 < c < c^{\ast}$, let 
\begin{equation}\label{eq:constforldv}
C(c,\Delta) = \left(1-e^{\log\Delta-I(c)}\right)^{-1/2}.
\end{equation}

Further, let $\lambda_2$ be the absolute value
of the second eigenvalue of the generator of the Glauber dynamics on
$G$, i.e.  $\lambda_2=\frac{1}{\tau_2}$. Let $f=f(\sigma)$ depend only
on the values of the configuration in $A$ and $g=g(\sigma)$ depend
only on the values of the configuration in $B$. If $\mu(f)=\mu(g)=0$,
then
\begin{equation}\label{eq:tgen}
\mu(fg)\leq \left(e^{-cd\lambda_2}+2C(c,\Delta)
\sqrt{\left|\partial_i A\right|e^{d\left(\log\Delta-I(c)\right)}}
\right)\|f\|_2\|g\|_2.
\end{equation}
In particular (by letting $c = e^{-\log\Delta-\gamma-2}$) for $\gamma
\geq 0$,
\begin{equation}\label{eq:tgen2}
\mu(fg)\leq \left(e^{-d\lambda_2\exp(-\log\Delta-\gamma-2)}+
4\sqrt{\exp(-(\gamma+1)d)\left|\partial_i A\right|}
\right)\|f\|_2\|g\|_2.
\end{equation}
\end{Proposition}

\begin{proof}[Proof of Theorem \ref{thm:general}]
  Note that $|\partial_i A|\leq|A|\leq\Delta^{r/2}$. Therefore, to prove
  (\ref{eq:l2_decay}) we use (\ref{eq:tgen2}) with B=$\{v:d(v,o)=r\}$,
  $d=r/2$ and $\gamma$ s.t. $e^\gamma>\Delta$.
\end{proof}

\begin{proof}[Proof of Proposition \ref{prop:gendecay}]
We use a  coupling argument.
Let $\mu$ be the Gibbs measure on $G$, and let $X_0$ be
 chosen according to $\mu$.
Let $X_t$ and $Y_t$ be defined as follows:
Set $Y_0=X_0$. For $t>0$,  let  $X_t$ and $Y_t$
evolve according to the dynamics with the following graphical
representation:
Each $v\in G$ has a Poisson clock. Assume the clock at $v$ rang at time $t$, and let $X_{t-}$ and $Y_{t-}$
be the configurations just before time $t$.
At time $t$ we do the following:
\begin{enumerate}
\item If $v\in B$ then $X_v$ updates according to the Gibbs measure,
  and $Y_v$ does not change.
\item
If $v\not\in B$ and $X_{t-}(w)=Y_{t-}(w)$ for every neighbor $w$ of $v$, then both $X$ and $Y$ update according
to the Gibbs measure so that $X_t(v)=Y_t(v)$.
\item
If $v\not\in B$ and there exists a neighbor $w$ of $v$ s.t. $X_{t-}(w)\neq Y_{t-}(w)$ then both $X$ and $Y$ update according
to the Gibbs measure, but this time independently of each other.
\end{enumerate}



For a vertex $v\in B$ we define $t_v$ to be the first time
the Poisson clock at $v$ rang.
For any $v \in G\setminus B$, we define
$t_v$ to be the first time the Poisson clock at $v$ rang after
$\min_{(w,v) \in E_G} t_w$.
Note that $X_t(v)=Y_t(v)$ at any time $t<t_v$,
and that
$t_v$ depends only on the Poisson clocks,
and is independent
of the initial configuration $X_0$.
We let $t_A = \min_{v \in A} t_v$.

\newcommand {\galp} {{\tilde{P}}}
\newcommand {\galq} {{\tilde{Q}}}

Let $\F_t$ denote the ($\sigma$-algebra of the) Poisson clocks at the vertices up to time $t$.
Let $(P^t f)(\sigma) = \E[f(X_t)|X_0=\sigma,\F_t]$ and
let $(Q^t f)(\sigma) = \E[f(Y_t)|X_0=\sigma,\F_t]$. Also,
let $(\galp^t f)(\sigma) = \E[f(X_t)|X_0=\sigma]$ and $(\galq^t f)(\sigma) = \E[f(Y_t)|X_0=\sigma]$.
Since for all $t$ the process $Y_t$ is at the stationary
distribution and $Y_t|_B=X_0|_B$ for all $t$, we get
\begin{equation} \label{eq:proj1}
\mu[g f] = \E[g(Y_t) f(Y_t)] =
\E[g(X_0) f(Y^t)] =
\E[g \galq^t f ].
\end{equation}

If $t<t_A$, then clearly $X_t=Y_t$ on $A$. Therefore, $\|(Q^t f-P^tf)\cdot
1_{t<t_A}\|_2 ^2=0$. On the other hand,
$\|Q^t f\|_2 \leq \|f\|_2 ~~ \F_t-\mbox{a.s.}$  and $\|P^t f\|_2 \leq \|f\|_2 ~~ \F_t-\mbox{a.s.}$
This is because the operators $f\to Q^t(f)$ and $f\to P^t(f)$  given $\F_t$ are Markov operators and hence
contractions. Therefore
\begin{eqnarray*}
\|\galq^t f-\galp^t f\|_2 ^2
&=&\E\left([\galq^t f(X_0)-\galp^t f(X_0)]^2\right)\\
& \leq &
    \E\left([Q^t f(X_0)-P^t f(X_0)]^2\right)\\
& = &
    \P(t\leq t_A)\int d\mu(\sigma)
\E\left([Q^t f(\sigma)-P^t f(\sigma)]^2 |{t \leq t_A} \right)\\
& + &
    \P(t > t_A)\int d\mu(\sigma)
\E\left([Q^t f(\sigma) -P^t f(\sigma) ]^2 |{t > t_A},X_0=\sigma\right)\\
& \leq &
    4 \P[t_A \leq t] \|f\|_2^2 \,
\end{eqnarray*}
where the first inequality is because $\galq^t f(X_0)-\galq^t f(X_0)$ is a conditional
expectation of $Q^t f(X_0) - P^t f(X_0)$,
and the second
inequality is because $\{t >  t_A\}$ is $\F_t$-measurable.
Therefore, by the Cauchy-Schwartz inequality,

\begin{equation} \label{eq:another_cs}
\E[(\galq^t f- \galp^t f)g]
 \leq 2 \sqrt{\P[t_A \leq t]} \|f\|_2\ \|g\|_2 \, .
\end{equation}
Since
\[
\E[g \galp^t f] \leq e^{-\lam_2t} \|f\|_2 \|g\|_2,
\]
We infer that from (\ref{eq:another_cs}) and (\ref{eq:proj1}) that
\begin{equation} \label{eq:coupling_and_mixing}
\mu[fg] \leq \left( e^{-\lam_2t} + 2 \sqrt{\P[t_A \leq t]} \right)
\|f\|_2 \|g\|_2 \, .
\end{equation}
It remains to bound the two terms in the right-hand side of
(\ref{eq:coupling_and_mixing}).

For $0<c<c^{\ast}$, we take $t = cd$.
We obtain that the first term is
$
e^{-cd\lam_2},
$ as desired.
It remains to bound $\P[t_A \leq t]$.
We note that $t_A \leq t$ only if there is some self-avoiding path
(sometimes referred to as ``path of disagreement")
between the $A$ and $B$
along which the discrepancy between the two
distributions has been conveyed in time less than $t$.

\ignore{
Note that there are at most $|A| (\Delta-1)^k$ such paths of length $k$ for all $k \geq r/2$.
We fix such a path $v_1,\ldots,v_k$ and bound the probability that this path was activated up to
time $t$.
This probability is clearly bounded by $\P[\Bin(t,\sGr^{-1}) \geq k]$ (think of ``success" as an
activation of the first non-active element of $v_1,\ldots,v_k$).
We let $c' > 0$ be a constant such that for all $m$ and $p$ and for
all $k \geq mp/c'$ one has the following tail
estimate:
\[
\P[\Bin(m,p) \geq k] \leq \Delta^{-2k}.
\]
(such a constant exists by standard large deviation estimates,
see, e.g., \cite[Corollary 2.4]{JLR}).
Thus,
\[
\P[\Bin(t,\sGr^{-1}) \geq k] \leq \Delta^{-2k}.
\]
So summing over all paths we obtain:
\begin{equation} \label{eq:exp_time}
\P[t_A \leq t] \leq |A| \sum_{k \geq r/2} (\Delta-1)^{k} \Delta^{-2k}.
\end{equation}
Thus both summands in  (\ref{eq:coupling_and_mixing}) decay
exponentially in $r$,
as claimed.
}

Time-reversing the process, this means that first-passage-percolation
with rate-$1$ exponential passage times starting at $A$ needs to
arrive at distance $d$ within time $cd$.  There are at most
$\left|\partial_i A\right|\Delta^{k}$ paths of length $k$ for the first-passage-percolation
for each $k\geq d$. Let $\tau(v,w)$ be the time needed to cross the edge
$(v,w)$. For each path $v_1,v_2,\ldots,v_k$,
\begin{eqnarray*}
\P\left(\tau(v_1,v_2)+\tau(v_2,v_3)+\ldots
  \tau(v_{k-1},v_k)<cd\right)
<e^{-kI(c)}
\end{eqnarray*}
where $I(c)=c-\log c-1$ is the large deviation rate function for the
exponential distribution.
Therefore,
\begin{eqnarray*}
\P(t_A\leq t)&\leq&\left|\partial_i A\right|\sum\limits_{k=d}^\infty\exp\left[k(-I(c')+\log\Delta)\right]\\
&\leq& C^2(c,\Delta)\left|\partial_i A\right|e^{d\left(\log\Delta-I(c)\right)}
\end{eqnarray*}
\end{proof}
Plugging this bound into (\ref{eq:coupling_and_mixing}), 
we obtain (\ref{eq:tgen}) as needed.

\section{Open Problems}
In this section we specify some relevant problems that are still open.
\begin{Problem}
  What is the relaxation time $\tau_2(n,b,b^{-1/2})$ of the Glauber
  dynamics of the Ising model on the $b$-ary tree of depth $n$ at the
  critical temperature $1-2\epsilon=\frac{1}{\sqrt{b}}$?
\end{Problem}
\noindent
Using the sum of spins as a test function, we learn that $\Omega(\log
n)$ is a lower bound for $\tau_2(n,b,b^{-1/2})$. We conjecture that
 the relaxation time is of order $\Theta(\log n)$.
A weaker conjecture is that
\[
\lim_{n\to\infty}\frac{\log(\tau_2(n,b,b^{-1/2}))}{n} = 0.
\]

\begin{Problem}
Fix $b$, and let
\[
\tau_2(\beta)=\lim_{n\to\infty}\frac{\log(\tau_2(n,b,\beta))}{n}.
\]
Theorem \ref{thm:tree} part \ref{treeupperbound} tells us that
$\tau_2(\beta)$ exists and is finite for all $\beta$. Show that
$\tau_2(\beta)$ is a monotone function of $\beta$.
This question is a special case of a more general monotonicity 
conjecture due to the fourth author, described in \cite{serban}.
See \cite{serban} where a monotonicity result is proven for the
Ising model on the cycle.
\end{Problem}

\begin{Problem}
For the Ising model (with free boundary conditions and no external
field) on a general graph of bounded degree, does the converse
of Theorem \ref{thm:general} hold, i.e., does uniform exponential decay
of point-to-set correlations imply a uniform spectral gap? \newline
(As pointed out by F. Martinelli (personal communication),
 the converse fails in certain lattices
if plus boundary conditions are allowed).
\end{Problem}

\begin{Problem}
Recall the general upper bound $n e^{(4 \cE(G) + 2 \Delta) \beta}$
on the relaxation time of Glauber dynamics
in terms of cut-width from  proposition \ref{prop:cut-width}.
For which graphs does a similar lower bound
of the form $\tau_2 \ge e^{c \cE(G) \beta}$
(for some constant $c>0$) hold at low temperature?

Such a lower bound is known to hold for boxes in a Euclidean
lattice, our results imply its validity for regular trees, and we
can also verify it for expander graphs. A  specific class of graphs
which could be considered here are the metric balls around a specific
vertex in an infinite graph $\Gamma$ that has critical probability
$p_c(\Gamma) <1$ for bond percolation.

\end{Problem}

\noindent{\bf Remark.} After the results presented here were
described in the extended abstract \cite{first},
striking further results on this topic were obtained by
F. Martinelli, A. Sinclair, and D. Weitz \cite{MSW}.
For the Ising model on regular trees, in the temperatures where
we show the Glauber dynamics has a uniform spectral gap, they show
it satisfies  a uniform log-Sobolev inequality; moreover, they
study in depth the effects of external fields and boundary conditions.

\noindent{\bf Acknowledgment.} We are grateful to David Aldous,
David Levin, Laurent Saloff-Coste and Peter Winkler for useful
discussions. We thank Dror Weitz for helpful comments on \cite{first}.


\begin{thebibliography}{99}
\bibitem{Io} Ioffe, D. (1996).  A note on the extremality of the disordered
state for the Ising model on the Bethe lattice.
{\em Lett.\  Math.\  Phys.} {\bf 37}, 137--143.
\end{thebibliography}
\end{document}